\documentclass[12pt,oneside,reqno]{amsart}
\usepackage{}
\usepackage{amssymb}
\usepackage{txfonts}
\usepackage{bbm}
\usepackage{cases}
\usepackage{amsmath}
\usepackage{graphicx}
\usepackage{mathrsfs}
\usepackage{stmaryrd}
\usepackage{amsfonts}
\usepackage{enumerate,amsmath,amssymb,amsthm}

\numberwithin{equation}{section}

\newcommand{\be}{\begin{eqnarray}}
\newcommand{\ee}{\end{eqnarray}}
\newcommand{\ce}{\begin{eqnarray*}}
\newcommand{\de}{\end{eqnarray*}}
\newtheorem{theorem}{Theorem}[section]
\newtheorem{lemma}[theorem]{Lemma}
\newtheorem{remark}[theorem]{Remark}
\newtheorem{definition}[theorem]{Definition}
\newtheorem{proposition}[theorem]{Proposition}
\newtheorem{example}[theorem]{Example}
\newtheorem{corollary}[theorem]{Corollary}

\def\e{{\mathrm{e}}}
\def\eps{\varepsilon}

\def\p{\partial}

\def\[{{\Big[}}
\def\]{{\Big]}}
\def\<{{\langle}}
\def\>{{\rangle}}
\def\({{\Big(}}
\def\){{\Big)}}

\def\bx{{\mathbf{x}}}

\def\dif{{\mathord{{\rm d}}}}

\def\no{\nonumber}
\def\={&\!\!=\!\!&}
\def\bt{\begin{theorem}}
\def\et{\end{theorem}}
\def\bl{\begin{lemma}}
\def\el{\end{lemma}}
\def\br{\begin{remark}}
\def\er{\end{remark}}

\def\bd{\begin{definition}}
\def\ed{\end{definition}}
\def\bp{\begin{proposition}}
\def\ep{\end{proposition}}
\def\bc{\begin{corollary}}
\def\ec{\end{corollary}}
\def\bx{\begin{example}}
\def\ex{\end{example}}

\def\cB{{\mathcal B}}

\def\cI{{\mathcal I}}
\def\cJ{{\mathcal J}}

\def\cT{{\mathcal T}}

\def\mE{{\mathbb E}}

\def\mK{{\mathbb K}}

\def\mN{{\mathbb N}}

\def\mR{{\mathbb R}}

\def\sB{{\mathscr B}}

\def\sD{{\mathscr D}}

\def\sK{{\mathscr K}}
\def\sL{{\mathscr L}}

\def\sQ{{\mathscr Q}}

\def\geq{\geqslant}
\def\leq{\leqslant}

\allowdisplaybreaks

\begin{document}

\title{Singular SDEs with critical non-local and non-symmetric L\'evy type generator}

\date{}

\author{Longjie Xie}

\dedicatory{
School of Mathematics and Statistics,
Jiangsu Normal University, Xuzhou, Jiangsu 221000, P.R.China\\
Emails: xlj.98@whu.edu.cn
}

\begin{abstract}
In this work, by using Levi's parametrix method we first construct the fundamental solution of the critical non-local operator perturbed by gradient. Then, we use the obtained estimates to prove the existence and uniqueness of strong solutions for stochastic differential equation driven by Markov process with irregular coefficients, whose generator is a non-local and non-symmetric L\'evy type operator.

  \bigskip

  \noindent{{\bf Keywords and Phrases:} Fundamental solution, non-local operator, pathwise uniqueness, singular SDEs}

\end{abstract}

\maketitle \rm

\section{Introduction and Main Results}

Consider the following critical non-local and non-symmetric operator perturbed by the gradient operator:
\begin{align}
\sL f(x):=\sL^{\kappa}f(x)+b(x)\cdot\nabla f(x),\quad\forall f\in C^\infty_0(\mR^d),  \label{LL}
\end{align}
where
$$
\sL^{\kappa}f(x):=\int_{\mR^d}\Big[f(x+z)-f(x)-1_{\{|z|\leq 1\}}z\cdot\nabla f(x)\Big]\frac{\kappa(x,z)}{|z|^{d+1}}\dif z.
$$
Here, $\kappa(x,z)$ is a measurable function on $\mR^d\times\mR^d$ satisfying
\begin{align}
0<\kappa_0\leq\kappa(x,z)\leq \kappa_1,\quad \kappa(x,z)=\kappa(x,-z),\quad\forall x,z\in\mR^d,    \label{kappa}
\end{align}
and for some $\beta\in(0,1)$,
\begin{align}
|\kappa(x,z)-\kappa(x',z)|\leq\kappa_2|x-x'|^\beta,\quad\forall x,x',z\in\mR^d,    \label{kappa2}
\end{align}
where $\kappa_0, \kappa_1, \kappa_2$ are positive constants. It is critical in the sense that the non-local operator $\sL^\kappa$ has the same order as the gradient operator $\nabla$. In particular, when $\kappa(x,z)\equiv c_da(x)$ is independent of $z$, we get $\sL^\kappa=a(x)\Delta^{\frac{1}{2}}$. Hence, the operator $\sL^\kappa$ can be seen as a generalization of the variable coefficient critical fractional Laplacian operator. The critical operator $\sL$ has particular interest in physics and mathematics (see \cite{Ca-Va,Si3} and references therein).
The symmetric in $z$ of $\kappa$ is a common assumption in the literature, see \cite{Ca-Si}. As a result, we can also write $\sL^\kappa$ as
\begin{align*}
\sL^\kappa f(x)&=\text{p.v.}\int_{\mR^d}\Big[f(x+z)-f(x)\Big]\frac{\kappa(x,z)}{|z|^{d+1}}\dif z\\
&=\frac{1}{2}\int_{\mR^d}\Big[f(x+z)+f(x-z)-2f(x)\Big]\frac{\kappa(x,z)}{|z|^{d+1}}\dif z,\quad\forall f\in C^\infty_0(\mR^d).
\end{align*}
The purpose of this paper is to study a jump type stochastic differential equation (SDE) with irregular coefficients as \eqref{sde1} below, which has the infinitesimal generator given exactly by (\ref{LL}) and whose driven noise is a family of pure jump Markov process, which can be even not L\'evy. This reflects the regularization effects of such kind of multiplicative noises on the deterministic system, see \cite{Fl}.

Our main tool to study the singular SDEs in this work will be the fundamental solution of the operator $\sL$.
Thus, we shall first construct the fundamental solution of $\sL$ by using the Levi's parametrix method. We remark that this part has independent interests. For $\beta\in (0,1)$, we introduce the usual H\"older space which is given by
$$
C^\beta_b(\mR^d):=\left\{f\in\sB(\mR^d): \|f\|_{C^\beta_b}:=\sup_{x\in\mR^d}|f(x)|
+\sup_{x\not= y\in\mR^d}\frac{|f(x)-f(y)|}{|x-y|^\beta}<\infty\right\}.
$$
The following is the first main result of this paper.

\bt\label{main1}
Assume that (\ref{kappa})-(\ref{kappa2}) hold and $b\in C^\beta_b(\mR^d)$ for some $\beta \in (0,1)$. Then there exists a unique transition density
function $p(t,x,y)$ on $\mR_+\times\mR^d\times\mR^d$ solving
\begin{align}
\p_tp(t,x,y)=\sL p(t,x,y),\quad x\neq y, \label{eq15}
\end{align}
and satisfying the following properties: for any $T>0$,
\begin{enumerate}[(i)]
\item (Upper bound) there exists a constant $c_1>0$ such that for all $t\in[0,T]$ and $x,y\in\mR^d$,
\begin{align}
p(t,x,y)\leq c_1t(|x-y|+t)^{-d-1}. \label{eq16}
\end{align}

\item (Gradient estimate) there is a constant $c_2>0$ so that for all $t\in[0,T]$ and $x,y\in\mR^d$,
\begin{align}
|\nabla_xp(t,x,y)|\leq c_2(|x-y|+t)^{-d-1}.\label{eq17}
\end{align}

\item (H\"older estimate) For any $\vartheta\in(0,\beta)$ and every $t\in[0,T]$, $x,x',y\in\mR^d$, it holds for some $c_3>0$ that
\begin{align}
|\nabla p(t,x,y)-\nabla p(t,x',y)|\leq c_3|x-x'|^{\vartheta} t^{-\vartheta}(|\tilde x-y|+t)^{-d-1},  \label{f2}
\end{align}
where $\tilde x$ is the one of the two points $x$ and $x'$ which is nearer to $y$.
\end{enumerate}
\et

The fundamental solutions (also called the heat kernel) of non-local operators have caused much attentions in the past decades, see \cite{Ch-Ki-So,Ch-Ku2,Ja,K-S-V} and references therein. Among all, we mention that in \cite{Bo-Ja0}, Bogdan and Jakubowski obtained the sharp two sided heat kernel estimates for the following perturbation of $\Delta^{\frac{\alpha}{2}}$
by gradient operator:
$$
\sL_{(\alpha)}f(x):=\Delta^{\frac{\alpha}{2}}f(x)+b(x)\cdot\nabla f(x),\ \ \alpha\in(1,2),
$$
where $b$ belongs to the Kato's class $\sK^{\alpha-1}_d$ defined as follows: for $\gamma>0$,
$$
\sK^{\gamma}_d:=\left\{f\in L^1_{loc}(\mR^d): \lim_{\eps\downarrow 0}\sup_{x\in\mR^d}
\int_{|x-y|\leq\eps}\frac{|f(y)|}{|x-y|^{d-\gamma}}\dif y=0\right\}.
$$
The reason of limiting $\alpha\in(1,2)$ lies in the fact that the heat kernel
of $\sL_{(\alpha)}$ is not comparable with that of $\Delta^{\frac{\alpha}{2}}$ for $\alpha\in(0,1)$ even when $b\equiv1$. In \cite{Xi-Zh}, Xie and Zhang studied
the fundamental solution of the critical case $a_t(x)\Delta^{1/2}+b_t(x)\cdot\nabla$ with coefficients in H\"older's space and obtained the sharp two sided estimates.
Recently, Chen and Zhang \cite{Ch-Zh} construct the fundamental solution
for the following nonlocal and non-symmetric operator:
$$
\sL^\kappa_{(\alpha)}f(x):=\int_{\mR^d}\Big[f(x+z)-f(x)-1_{\{|z|\leq 1\}}z\cdot\nabla f(x)\Big]\frac{\kappa(x,z)}{|z|^{d+\alpha}}\dif z,
$$
where $\kappa$ satisfies (\ref{kappa})-(\ref{kappa2}) and $\alpha\in(0,2)$.

Our results can be seen as a generalization of \cite{Xi-Zh} to the more general non-local operators $\sL^\kappa$, as well as a generalization of \cite{Ch-Zh} to the drift perturbation in the critical case. We point out that the estimate (\ref{f2}) seems to be new in all the works mentioned above, which means that the function $p(t,x,y)$ has $``1+\vartheta"$- order regularity in $x$ with $\vartheta<\beta$. This regularity estimate is certainly delicate than (\ref{eq17}) and the proof is much more involved, as we shall see. What is more, it turns out to be of critical important below for us to study the singular SDEs by using the heat kernel estimates.

\vspace{3mm}
As mentioned above, we are led to the study of this non-local operator $\sL$ by the consideration of a kind of SDEs on $\mR^d$ driven by Markov process. To specify the SDE that we are going to study, denote by $m$ the Lebesgue measure, and let $N$ be a Poisson random measure on $\mR^d\times  [0,\infty)\times [0,\infty)$ with mean measure $\nu\times m\times m$, where $\nu$ is a L\'evy measure of Cauchy-type satisfying
\begin{align}
\nu(\dif z)=\frac{\bar\kappa(z)}{|z|^{d+1}}\dif z, \quad \bar\kappa(z)=\bar\kappa(-z),\quad \bar\kappa_0\leq \bar\kappa(z)\leq \bar\kappa_1,\label{nu}
\end{align}
here, $\bar\kappa(z)$ is a measurable function on $\mR^d$ and $\bar\kappa_0, \bar\kappa_1$ are two positive constants.
Set for $A\in \sB\big(\mR^d\times [0,\infty)\times  [0,\infty)\big)$,
$$
\tilde N(A):=N(A)-\nu\times m\times m(A).
$$
Consider the following SDE:
\begin{align}
\dif X_t=\int_0^\infty\!\!\!\!\int_{|z|\leq 1}\!&1_{[0,\sigma(X_{t-},z)]}(r)z\tilde N(\dif z\times\dif r\times\dif t)\no\\
&+\int_0^\infty\!\!\!\!\int_{|z|>1}\!1_{[0,\sigma(X_{t-},z)]}(r)z N(\dif z\times\dif r\times\dif t)+b(X_t)\dif t,\quad X_0=x\in\mR^d. \label{sde1}
\end{align}
An application of It\^o's formula shows that the generator is
\begin{align*}
\hat\sL f(x)
=\int_{\mR^d}\Big[f(x+z)-f(x)-1_{\{|z|\leq 1\}}z\cdot\nabla f(x)\Big]\sigma(x,z)\nu(\dif z)+b(x)\cdot\nabla f(x).
\end{align*}
If we let
$$
\hat\kappa(x,z):=\sigma(x,z)\bar\kappa(z),
$$
then, we can get
\begin{align}
\hat\sL f(x)=\sL^{\hat\kappa}f(x)+b(x)\cdot\nabla f(x).  \label{LLLL}
\end{align}
Hence, the generator of the above SDE is given exactly by $\sL$ as in (\ref{LL}). This makes (\ref{sde1}) more interesting and is worthy of study.

Under the conditions that $b$ is bounded and global Lipschitz continuous, $\sigma$ satisfies
\begin{align}
\int_{\mR^d}|\sigma(x,z)-\sigma(y,z)|\cdot|z|\nu(\dif z)\leq C_1|x-y|, \label{kur}
\end{align}
and with some other assumptions, it was showed by Kurtz \cite[Theorem 3.1]{Kurz2} that \eqref{sde1} has a unique strong solution, see also \cite{Ku-Po}. We shall study the pathwise uniqueness of strong solutions to SDE (\ref{sde1}) with irregular coefficients.
Using the conclusions obtained in Theorem \ref{main1}, we have the following result.

\bt\label{main2}
Let $\sigma(x,z)$ satisfies (\ref{kappa})-(\ref{kappa2}) with $\beta>\tfrac{1}{2}$ and the L\'evy measure $\nu$ satisfies (\ref{nu}). Suppose also that:

\begin{enumerate}
\item [{\bf (H$\sigma$)}] There exists a Kato function $h\in \mK_d^1$ (see Definition \ref{Def1} below) such that for almost all $x,y\in \mR^d$,
\begin{align}
\int_{\mR^d}|\sigma(x,z)-\sigma(y,z)|(|z|\wedge1)\nu(\dif z)\leq |x-y|\Big(h(x)+h(y)\Big).\label{a1}
\end{align}

\item [{\bf (Hb)}] For some $\theta\in (\tfrac{1}{2},1)$, the drift $b\in C^\theta_b(\mR^d)$.
\end{enumerate}
Then, for each $x\in\mR^d$, there exists a unique strong solution $X_t(x)$ to SDE (\ref{sde1}). Moreover, $X_t(x)$ admits a density function $p(t,x,y)$ which enjoys all the properties stated in the conclusions of Theorem \ref{main1}.
\et

Notice that the drift coefficient is singular enough so that the deterministic ordinary equation of (\ref{sde1}) with $\sigma\equiv0$ is not well-posed. The fact that noises may produce regularization effects which make ill-posed deterministic systems well-posed has attracted a lot of attentions in the past decades. A remarkable result of N. V. Krylov and M. R\"ockner \cite{Kr-Ro} shows that under the condition $b\in L^q_{loc}\big(\mR_+;L^p(\mR^d)\big)$ with
\begin{align}
d/p+2/q<1,\label{pq}
\end{align}
there exists a unique strong solution for every $x\in\mR^d$ to the following SDE:
$$
\dif X_t=\dif W_t+b(t,X_t)\dif t,\quad X_0=x\in\mR^d.
$$
Later on, this was extended by Zhang \cite{Zh3} to the multiplicative noise
\begin{equation}  \label{e:ZhEqn}
\dif X_t=\sigma(t,X_t)\dif W_t+b(t,X_t)\dif t,\quad X_0=x\in\mR^d
\end{equation}
under the assumption that $\sigma$ is uniformly continuous in $x$, bounded and uniformly elliptic and
\begin{align}
\nabla\sigma\in L^q_{loc}\big(\mR_+;L^p(\mR^d)\big) \label{multi}
\end{align}
with $p,q$ satisfy (\ref{pq}), where $\nabla$ donets the weak detivative of $\sigma$ with respect to $x$. See also \cite{Fe-Fl-2,W-Z,XZ}. The situation for SDEs with pure jump L\'evy noises is more delicate. Let $(L_t)_{t\geq 0}$ be a symmetric $\alpha$-stable process with non-degenerate spectral measure and consider the following SDE:
\begin{align}
\dif X_t=\dif L_t+b(X_t)\dif t,\quad X_0=x\in\mR^d. \label{levy}
\end{align}
When $\alpha\geq 1$ and $b$ is $\beta$-H\"older continuous with
$$
\beta>1-\frac{\alpha}{2},
$$
it was proved by Priola \cite{Pri} that there exists a unique strong solution $X_t(x)$ to SDE (\ref{levy}) for each $x\in\mR^d$. Zhang \cite{Zh00} extended this result when $\alpha>1$ and allowing $b$ in some fractional Sobolev space. See also \cite{B-B-C,Ch-So-Zh,Pri2} for related results. Recently, \cite{Xi-Xu} considered the same SDEs as (\ref{sde1}) with the L\'evy measure given by
$\nu(\dif z)=\bar\kappa(z)|z|^{-d-\alpha}\dif z,$
where $\bar\kappa$ satisfies (\ref{nu}) and $\alpha\in(1,2)$. Thus, we fill the gap in the critical case $\alpha=1$ in this paper.

\vspace{4mm}
Compared with \cite{Kr-Ro,Pri,Zh3,Zh00}, we shall use directly the estimates of fundamental solution in the whole procedure. It seems the first time that Kato functions which are commonly used in the study of heat kernel estimates are bringed to the study of strong solutions for singular SDEs. Our approach can also be adapted to SDEs driven by multiplicative Brownian motion. Another advantage of our method is that, as an consequence of (\ref{eq17}) and Theorem \ref{main2}, we can derive the following estimate of the semigroup corresponding to $X_t$.

\bc
Let $X_t$ be the unique strong solution to SDE (\ref{sde1}) and $\cT_t$ be the corresponding semigroup. Then, $X_t$ is strong Feller and
\begin{align}
|\nabla \cT_tf(x)|\leq Ct^{-\frac{1}{\alpha}}\|f\|_{\infty},\quad\forall f\in\cB_b(\mR^d),  \label{88}
\end{align}
where $C$ is a positive constant.
\ec

\br
In general, for SDE (\ref{sde1}), it is not easy to deduce the strong Feller property as well as the derivative formula for the corresponding semigroup even if the coefficients $\sigma$ and $b$ belong to $C_b^\infty(\mR^d)$. The trouble is caused by the term $1_{[0,\sigma(x,z)]}$, which is not differentiable even though $\sigma$ is smooth. Here, we easily get the estimate (\ref{88}) through the fundamental solution estimates.
\er

Let us specify the main difficulties of the proof.
As is well known now, the basic idea of the proof for the pathwise uniqueness of singular SDEs is based on the the Zvonkin's transformation, which require suitable analytic regularity results of certain elliptic equations corresponding to the generators of the strong solutions. Notice that \cite{Pri,Pri2,Zh00} are all restricted to the additive L\'evy noise. In this circumstances, one only needs to deal with the symmetric operator
$\sL_0$ defined by
$$
\sL_0f(x):=\int_{\mR^d}\Big[f(x+z)-f(x)-1_{\{|z|\leq 1\}}z\cdot\nabla f(x)\Big]\nu(\dif z),\quad\forall f\in C_0^\infty(\mR^d),
$$
which is the generator of $L_t$. The analysis in \cite{Pri,Pri2,Zh00} relies heavily on
the symmetric property of $\sL_0$ and the $C^2$-smoothing property of its semigroup.
However, the operator $\sL^\kappa$ in our paper is not symmetric, critical and more important, its semigroup has only $C^{1+\vartheta}$ regularity as indicated by (\ref{f2}). Therefore, we need to use more delicate analysis to fit our less regularity property into the frame of Zvonkin's transformation.
Another difficult comes from the new extra term $1_{[0,\sigma(X_{s-},z)]}(r)$. We need to use a trick of $L_1$-estimate by Kurtz \cite{Kurz2} rather than the $L_2$-estimate as usual when proving our main theorem. Some new challenges appear when dealing with the $L_1$-estimate and the irregular coefficients, see also \cite{Xi-Xu}.

\vspace{4mm}
Last but not least, it is clear that the assumption (\ref{a1}) is a generalization of (\ref{kur}). Here, we would like to give the following important comment.

\br
In view of (\ref{incl}) below, we can take
$$
\sigma(x,z)=K(z)+\tilde\sigma(x)|z|^{\gamma} \ \ {\rm for} \ \  |z| \le 1, \ \ \ \ \ \sigma(x,z)=K(z)+\tilde\sigma(x) \ \ {\rm for} \ \  |z|>1,
$$
with $0<K_1\leq K(z)\leq K_2$, $\gamma>0$ and $\nabla\tilde\sigma\in L^q_{loc}(\mR^d)$ with $q>d$, where $\nabla$ denotes the weak derivative.
The interesting thing is that for SDE \eqref{e:ZhEqn}, if $\sigma$ is independent of the time variable $t$, one has to assume that $\nabla\sigma\in L^q(\mR^d)$ with $q>d$, see (\ref{pq}) and (\ref{multi}). So one may guess reasonably that for SDEs driven by multiplicative $\alpha$-stable noise, one has to assume that the diffusion coefficients satisfies $\nabla \sigma\in L^{q}(\mR^d)$ with $q>2d/\alpha$. Here, we only need the index $q>d$, the point is that $\sigma$ appears in the indicator function $1_{[0,\sigma(X_{s-},z)]}(r)$.
\er

This paper proceed as follows: In Section 2 we construct the fundamental solution for the operator $\sL$ and give the proof of Theorem \ref{main1}. Meanwhile, we study the smoothing properties of the corresponding semigroup, which will play an essential role. In Section 3, we prove our main result Theorem \ref{main2}.
Throughout this paper, we use the following convention: $C$ with or without subscripts will denote a positive constant, whose value may change in different places, and whose dependence on parameters can be traced from calculations. We write $f(x)\preceq g(x)$ to mean that there exists a constant
$C_0>0$ such that $f(x)\leq C_0 g(x)$; and $f(x)\asymp g(x)$ to mean
that there exist $C_1,C_2>0$ such that $C_1 g(x)\leq f(x)\leq C_2 g(x)$.

\section{Fundamental solution of operator $\sL$}

\subsection{Preliminaries}
To shorten the notation, we set for $\gamma, \beta\in \mR$,
$$
\varrho_\gamma^\beta(t,x):=t^{\gamma}\big(|x|^{\beta}\wedge 1\big)\big(|x|+t\big)^{-d-1}.
$$
Let $\cB(\gamma,\beta)$ be the usual Beta function defined by
$$
\cB(\gamma,\beta):=\int^1_0(1-s)^{\gamma-1}s^{\beta-1}\dif s,\ \
\gamma,\beta>0.
$$
The following result which is called the $3$P-inequality was proved in \cite[Lemma 2.1]{Xi-Zh}.
\bl\label{3p}
For $\beta_1,\beta_2\in [0,1]$ and $\gamma_1,\gamma_2\in\mR$, there exists a constant $C_d>0$ such that for all
$0\leq s, t<\infty$ and $x,y\in\mR^d$,
\begin{align}
\int_{\mR^d}\varrho^{\beta_1}_{\gamma_1}(t,x-z)\varrho^{\beta_2}_{\gamma_2}(s,z-y)\dif z
&\leq C_d\Big\{t^{\gamma_1+\beta_1+\beta_2-1}s^{\gamma_2}\varrho^0_0(t+s,x-y)\no\\
&\quad+t^{\gamma_1+\beta_1-1}s^{\gamma_2}\varrho^{\beta_2}_0(t+s,x-y)\no\\
&\quad+t^{\gamma_1}s^{\gamma_2+\beta_1+\beta_2-1}\varrho^0_0(t+s,x-y)\no\\
&\quad+t^{\gamma_1}s^{\gamma_2+\beta_2-1}\varrho^{\beta_1}_0(t+s,x-y)\Big\},\label{EU7}
\end{align}
and if $\gamma_1>-\beta_1$, $\gamma_2>-\beta_2$, we have
\begin{align*}
\int_0^t\!\!\!\int_{\mR^d}\varrho^{\beta_1}_{\gamma_1}(t-s,x,z)\varrho^{\beta_2}_{\gamma_2}(s,z,y)\dif z\dif s
&\leq C_d\Big\{\varrho^0_{\gamma_1+\gamma_2+\beta_1+\beta_2}(t,x,y)\cB(\gamma_1+\beta_1+\beta_2,1+\gamma_2)\\
&\quad+\varrho^{\beta_2}_{\gamma_1+\gamma_2+\beta_1}(t,x,y)\cB(\gamma_1+\beta_1,1+\gamma_2)\\
&\quad+\varrho^0_{\gamma_1+\gamma_2+\beta_1+\beta_2}(t,x,y)\cB(\gamma_2+\beta_1+\beta_2,1+\gamma_1)\\
&\quad+\varrho^{\beta_1}_{\gamma_1+\gamma_2+\beta_2}(t,x,y)\cB(\gamma_2+\beta_2,1+\gamma_1)\Big\}.
\end{align*}
\el

Let us first recall some facts about the heat kernel of the following non-local symmetric operator (with a little abuse of notation, we still denote it by $\sL^\kappa$):
\begin{align*}
\sL^\kappa f(x):=\text{p.v.}\int_{\mR^d}\Big[f(x+z)-f(x)\Big]\kappa(z)|z|^{-d-1}\dif z,\quad \forall f\in C_0^\infty(\mR^d).
\end{align*}
Here, the function $\kappa$ is independent of $x$ and satisfies
\begin{align}
\kappa(z)=\kappa(-z),\quad 0<\kappa_0\leq \kappa(z)\leq \kappa_1,\quad \forall z\in\mR^d.    \label{kap1}
\end{align}
It is known that there exists a symmetric $\alpha$-stable like process on $\mR^d$ corresponding to $\sL^{\kappa}$. Let $Z^\kappa(t,x)$ be the heat kernel of operator $\sL^{\kappa}$, i.e.,
\begin{align*}
\p_tZ^\kappa(t,x)=\sL^\kappa Z^\kappa(t,x),\quad \lim_{t\downarrow0}Z^\kappa(t,x)=\delta_0(x),
\end{align*}
where $\delta_0(x)$ is the Dirac function. Then, it follows from \cite[Theorem 1.1]{Ch-Ku} that for some constant $C_0>1$,
\begin{align}
C_0^{-1}\varrho_1^0(t,x)\leq Z^\kappa(t,x)\leq C_0\varrho_1^0(t,x),\quad\forall t\geq 0,\,\,x\in\mR^d. \label{p0}
\end{align}
Moreover, if we set
\begin{align}
\hat\kappa(z):=\kappa(z)-\frac{\kappa_0}{2},  \label{pp1}
\end{align}
by the construction of L\'evy process, we can write (see also \cite[(2.23)]{Ch-Zh})
\begin{align}
Z^\kappa(t,x)=\int_{\mR^d}\rho(\tfrac{\kappa_0}{2}t,x-z)Z^{\hat\kappa}(t,z)\dif z, \label{kk}
\end{align}
where $\rho$ is the heat kernel of $\Delta^{\frac{1}{2}}$ given by
$$
\rho(t,x)=\pi^{-\frac{d+1}{2}}\Gamma(\tfrac{d+1}{2})(|x|^2+t^2)^{-\frac{d+1}{2}}t,
$$
and $\Gamma$ is the usual Gamma function. This is also called the Poisson kernel.

Below, for a function $f$ on $\mR_+\times\mR^d$, we shall simply write
$$
\delta_f(t,x;z):=f(t,x+z)+f(t,x-z)-f(t,x).
$$
By \cite[Lemma 2.3, Theorem 2.4]{Ch-Zh}, it holds for all $t\geq 0, x\in\mR^d$ that there exist positive constants $C_1,C_2$ such that
\begin{align}
|\nabla Z^\kappa(t,x)|\leq C_1\varrho_0^0(t,x) \label{p1}
\end{align}
and
\begin{align}
\int_{\mR^d}\big|\delta_{Z^\kappa}(t,x;z)\big|\cdot|z|^{-d-1}\dif z\leq C_2\varrho_{0}^0(t,x), \label{p2}
\end{align}
Moreover, we also have the following H\"older estimates:
for any $\vartheta\in(0,1]$, $t>0$ and all $x,x'\in\mR^d$, there exist $C_3,C_4>0$ such that
\begin{align}
|\nabla Z^\kappa(t,x)-\nabla Z^\kappa(t,x')|\leq C_3|x-x'|^\vartheta\varrho_{-\vartheta}^0(t,\tilde x) \label{zh}
\end{align}
and
\begin{align}
\int_{\mR^d}\big|\delta_{Z^\kappa}(t,x;z)-\delta_{Z^\kappa}(t,x';z)\big|\cdot|z|^{-d-1}\dif z\leq C_4|x-x'|^\vartheta\varrho_{-\vartheta}^0(t,\tilde x), \label{p3}
\end{align}
where $\tilde x$ is the one of the two points $x$ and $x'$ which is nearer to zero point. In fact, (\ref{zh}) is shown by \cite[Lemma 2.3]{Xi-Xu}. As for (\ref{p3}), if $|x-x'|>t$, then it is easy to see by (\ref{p2}) that
\begin{align*}
\int_{\mR^d}\big|\delta_{Z^\kappa}(t,x;z)-\delta_{Z^\kappa}(t,x';z)\big|\cdot|z|^{-d-1}\dif z\leq C_4\Big(\varrho_0^0(t,x)+\varrho_0^0(t,x')\Big)\leq C_4|x-x'|^\vartheta\varrho_{-\vartheta}^0(t,\tilde x).
\end{align*}
In the case $|x-x'|\leq t$, we use \cite[Theorem 2.4]{Ch-Zh} to deduce
\begin{align*}
\int_{\mR^d}\big|\delta_{Z^\kappa}(t,x;z)-\delta_{Z^\kappa}(t,x';z)\big|\cdot|z|^{-d-1}\dif z\leq C_4|x-x'|\Big(\varrho_{-1}^0(t,x)+\varrho_{-1}^0(t,x')\Big)\leq C_4|x-x'|^\vartheta\varrho_{-\vartheta}^0(t,\tilde x),
\end{align*}
thus (\ref{p3}) is true. Let $\kappa$ and $\tilde\kappa$ be two functions on $\mR^d$ satisfying (\ref{kap1}), we shall also need the following continuous dependence of the heat kernel with respect to the kernel function $\kappa$:
\begin{align}
|Z^\kappa(t,x)-Z^{\tilde\kappa}(t,x)|  \leq C_5\|\kappa-\tilde\kappa\|_{\infty}\Big(\varrho_1^0+\varrho_{1-\gamma}^\gamma\Big)(t,x),  \label{con1}
\end{align}
and
\begin{align}
|\nabla Z^\kappa(t,x)-\nabla Z^{\tilde\kappa}(t,x)|  \leq C_6\|\kappa-\tilde\kappa\|_{\infty}\Big(\varrho_0^0+\varrho_{-\gamma}^\gamma\Big)(t,x),  \label{con}
\end{align}
where $C_5,C_6>0$ are constant and $\gamma\in(0,1)$, see \cite[Theorem 2.5]{Ch-Zh}.

\subsection{Construction and estimates of the fundamental solution}
Now, we consider the operator $\sL$ in (\ref{LL}), which can be seen as $\sL^\kappa$ perturbed by the gradient term. In order to reflect the dependence of $\kappa$ with respect to $x$, we also write
$$
\sL^{\kappa(x)}f(x)=\sL^\kappa f(x).
$$
Notice that the operator $\sL^\kappa$ has the same order with $\nabla$. Hence, the usual perturbation method to construct the heat kernel is not applicable. As in \cite{Xi-Zh}, we shall use the Levi's parametrix method. Fix $y\in\mR^d$, consider the freezing operator
$$
\sL^{\kappa(y)}f(x):=\text{p.v.}\int_{\mR^d}\Big[f(x+z)-f(x)\Big]\kappa(y,z)|z|^{-d-1}\dif z.
$$
Let $Z_y(t,x):=Z^{\kappa(y)}(t,x)$ be the heat kernel of operator $\sL^{\kappa(y)}$, i.e.,
\begin{align}
\p_tZ_y(t,x)=\sL^{\kappa(y)} Z_y(t,x),\quad \lim_{t\downarrow0}Z_y(t,x)=\delta_0(x).   \label{pp2}
\end{align}
For a bounded and measurable function $b$, we define
$$
p_0(t,x,y):=Z_y\big(t,x-y+b(y)t\big).
$$
Then, one can check by (\ref{pp2}) that
\begin{align}
\p_tp_0(t,x,y)=\sL^{\kappa(y)}p_0(t,x,y)+b(y)\cdot\nabla p_0(t,x,y),\quad \lim_{t\downarrow0}p_0(t,x,y)=\delta_y(x).   \label{pp3}
\end{align}
Meanwhile, we have the following important estimates, which will be used below.

\bl
Under (\ref{kappa}), there exist constants $C_0>1, C_1>0$ such that for every $t\geq 0$ and $x,y\in\mR^d$,
\begin{align}
C_0^{-1}\varrho_1^0(t,x-y)\leq p_0(t,x,y)\leq C_0\varrho_1^0(t,x-y)   \label{p00}
\end{align}
and
\begin{align}
|\nabla p_0(t,x,y)|+\int_{\mR^d}\big|\delta_{p_0}(t,x,y;z)\big|\cdot|z|^{-d-1}\dif z\leq C_1\varrho_{0}^0(t,x-y). \label{frp0}
\end{align}
Moreover, for any $\vartheta\in(0,1]$, there exist $C_2,C_3>0$ such that
\begin{align}
|\nabla p_0(t,x,y)-\nabla p_0(t,x',y)|\leq C_2|x-x'|^\vartheta\varrho_{-\vartheta}^0(t,\tilde x-y)    \label{p02}
\end{align}
and
\begin{align}
\int_{\mR^d}\big|\delta_{p_0}(t,x,y;z)-\delta_{p_0}(t,x',y;z)\big|\cdot|z|^{-d-1}\dif z\leq C_3|x-x'|^\vartheta\varrho_{-\vartheta}^0(t,\tilde x-y), \label{p03}
\end{align}
where $\tilde x$ is the one of the two points $x$ and $x'$ which is nearer to $y$.
\el

\begin{proof}
Since $\kappa$ is uniformly bounded, it follows by (\ref{p0}) and the definition of $p_0$ that for some constant $C_0$ independent of $y$,
\begin{align*}
p_0(t,x,y)\asymp \varrho_1^0\big(t,x-y+b(y)t\big).
\end{align*}
Noticing that $b$ is bounded,
$$
|x-y+b(y)t|+t\asymp|x-y|+t,
$$
we get (\ref{p00}). Similarly, by (\ref{p1})-(\ref{p3}) we can get (\ref{frp0})-(\ref{p03}).
\end{proof}

We also prepare the following important estimates for latter use.

\bl
Assume that (\ref{kappa}) holds and $b\in C_b^\beta(\mR^d)$ for some $\beta\in(0,1)$. We have for all $t\geq 0$ and $x\in\mR^d$,
\begin{align}
\bigg|\int_{\mR^d}\nabla p_0(t,x,y)\dif y\bigg|\leq C_{d}t^{\beta-1},\label{00}
\end{align}
and for any $\vartheta\in(0,1)$ and $x, x'\in\mR^d$,
\begin{align}
\bigg|\int_{\mR^d}\Big[\nabla p_0(t,x,y)-\nabla p_0(t,x',y)\Big]\dif y\bigg|\leq C_{d,\vartheta}|x-x'|^\vartheta t^{\beta-\vartheta-1},\label{000}
\end{align}
where $C_{d}, C_{d,\vartheta}$ are positive constants.
\el

\begin{proof}
Since $Z_y(t,x)$ is the heat kernel of the operator $\sL^{\kappa(y)}$, we have
$$
\int_{\mR^d}Z_y(t,x)\dif x=1,\quad\forall y\in\mR^d.
$$
As a result, we can also get
$$
\int_{\mR^d}Z_\xi\big(t,x-y+b(\xi)t\big)\dif y=1,\quad\forall x,\xi\in\mR^d.
$$
In view of (\ref{con}) and using (\ref{p02}) with $\vartheta=1$, we find that for any $\gamma\in(0,1)$,
\begin{align}
&\big|\nabla Z_y\big(t,x-y+b(y)t\big)-\nabla Z_\xi\big(t,x-y+b(\xi)t\big)\big|\no\\
&\leq \big|\nabla Z_y\big(t,x-y+b(y)t\big)-\nabla Z_\xi\big(t,x-y+b(y)t\big)\big|\no\\
&\quad+\big|\nabla Z_\xi\big(t,x-y+b(y)t\big)-\nabla Z_\xi\big(t,x-y+b(\xi)t\big)\big|\no\\
&\preceq \big(|\xi-y|^{\beta}\wedge1\big)\Big(\varrho_0^0+\varrho_{-\gamma}^\gamma\Big)(t,x-y)+\big(|\xi-y|^{\beta}\wedge1\big)t\cdot\varrho_{-1}^0(t,x-y)\no\\
&\preceq\big(|\xi-y|^{\beta}\wedge1\big)\Big(\varrho_0^0+\varrho_{-\gamma}^\gamma\Big)(t,x-y).\label{invo}
\end{align}
Hence, we can deduce
\begin{align*}
\int_{\mR^d}\nabla p_0(t,x,y)\dif y&=\int_{\mR^d}\Big[\nabla p_0(t,x,y)-\nabla Z_\xi\big(t,x-y+b(\xi)t\big)\Big]\Big|_{\xi=x}\dif y\\
&\leq \int_{\mR^d}\Big(\varrho_0^\beta+\varrho_{-\gamma}^{\gamma+\beta}\Big)(t,x-y)\dif y\preceq t^{\beta-1},
\end{align*}
which gives (\ref{00}). The estimate (\ref{000}) is more involved. Let $\hat Z_y(t,x)$ be the heat kernel of operator $\sL^{\hat\kappa(y)}$, where $\hat\kappa(y)$ is defined as in (\ref{pp1}). By (\ref{kk}), we can write
\begin{align*}
p_0(t,x,y)&=\int_{\mR^d}\rho\big(\tfrac{\kappa_0}{2}t,x-y+b(y)t-z\big)\hat Z_y(t,z)\dif z\\
&=\int_{\mR^d}\rho\big(\tfrac{\kappa_0}{2}t,x-z\big)\hat Z_y\big(t,z-y+b(y)t\big)\dif z.
\end{align*}
For simplicity, we set
$$
\zeta_{\nabla \rho}(t;x,x';z):=\nabla \rho\big(\tfrac{\kappa_0}{2}t,x-z\big)-\nabla \rho\big(\tfrac{\kappa_0}{2}t,x'-z\big).
$$
and let $\tilde x$ be the one of the two points $x$ and $x'$ which is nearer to $z$. Then, we know as in (\ref{zh}) that for any $\vartheta\in(0,1)$,
$$
|\zeta_{\nabla \rho}(t;x,x';z)|\leq C_\vartheta|x-x'|^\vartheta\varrho_{-\vartheta}^0\big(t,\tilde x-z\big).
$$
We may argue as above to deduce that
\begin{align*}
\cJ&:=\int_{\mR^d}\Big[\nabla p_0(t,x,y)-\nabla p_0(t,x',y)\Big]\dif y\\
&=\int_{\mR^d}\!\int_{\mR^d}\zeta_{\nabla \rho}(t;x,x';z)\hat Z_y\big(t,z-y+b(y)t\big)\dif z\dif y\\
&=\int_{\mR^d}\!\int_{\mR^d}\zeta_{\nabla \rho}(t;x,x';z)\big[\hat Z_y\big(t,z-y+b(y)t\big)-\hat Z_\xi\big(t,z-y+b(\xi)t\big)\big]\Big|_{\xi=\tilde x}\dif z\dif y.
\end{align*}
As in (\ref{invo}), we use (\ref{con1}), (\ref{p00}) and (\ref{frp0}) to deduce that for any $0<\gamma<1$, there exists a $C_{\gamma}$ such that
\begin{align*}
\big[\hat Z_y\big(t,z-y+b(y)t\big)-\hat Z_\xi\big(t,z-y+b(\xi)t\big)\big]\Big|_{\xi=\tilde x}\leq C_{\gamma}\big(|\tilde x-y|^{\beta}\wedge1\big)\Big(\varrho^{0}_{1}+\varrho^{\gamma}_{1-\gamma}\Big)(t,z-y),
\end{align*}
which yields by (\ref{EU7}) that
\begin{align*}
\cJ&\leq C_\vartheta|x-x'|^\vartheta\int_{\mR^d}\!\int_{\mR^d}\varrho^0_{-\vartheta}(t,\tilde x-z)\Big(\varrho^{0}_{1}+\varrho^{\gamma}_{1-\gamma}\Big)(t,z-y)\dif z\cdot\big(|\tilde x-y|^{\beta}\wedge1\big)\dif y\\
&\leq C_\vartheta|x-x'|^\vartheta\int_{\mR^d}\Big(\varrho^{\beta}_{-\vartheta}+\varrho^{\gamma+\beta}_{-\vartheta-\gamma}\Big)\big(t,\tilde x-y\big)\dif y\leq C_\vartheta|x-x'|^\vartheta t^{\beta-\vartheta-1}.
\end{align*}
The proof is finished.
\end{proof}

By Levi's parametrix method, we construct the fundamental solution of $\sL$ by the following formula:
\begin{align}
p(t,x,y)=p_0(t,x,y)+\int_0^t\!\!\!\int_{\mR^d}p_0(t-s,x,z)q(s,z,y)\dif z\dif s, \label{defp}
\end{align}
where $q$ satisfies the following integral equation:
\begin{align}
q(t,x,y)=q_0(t,x,y)+\int_0^t\!\!\!\int_{\mR^d}q_0(t-s,x,z)q(s,z,y)\dif z\dif s, \label{q}
\end{align}
and
$$
q_0(t,x,y):=\Big(\sL^{\kappa(x)}-\sL^{\kappa(y)}\Big)p_0(t,x,y)+\Big(b(x)-b(y)\Big)\cdot\nabla p_0(t,x,y).
$$

\br
The point is that we should freeze simultaneously both the diffusion coefficient and the drift coefficient at the given point $y$, see the definition of $p_0$.
\er

Formally, we have by (\ref{pp3}) that
\begin{align}
\p_tp(t,x,y)&=\p_tp_0(t,x,y)+q(t,x,y)+\int_0^t\!\!\!\int_{\mR^d}\p_tp_0(t-s,x,z)q(s,z,y)\dif z\dif s \no\\
&=\sL^{\kappa(x)} p_0(t,x,y)+b(x)\cdot\nabla p_0(t,x,y)\no\\
&\qquad\qquad\quad+\int_0^t\!\!\!\int_{\mR^d}\Big(\sL^{\kappa(x)} p_0(t-s,x,z)+b(x)\cdot\nabla p_0(t-s,x,z)\Big)q(s,z,y)\dif z\dif s \no\\
&=\sL^{\kappa(x)} p(t,x,y)+b(x)\cdot\nabla p(t,x,y). \label{form}
\end{align}
Thus, the main tasks are to solve the equation (\ref{defp}), and to make the above computations rigorous.

Below, we shall work on the time interval $[0,1]$, and always assume
(\ref{kappa})-(\ref{kappa2}) hold and $b\in C^\beta_b(\mR^d)$.  The general case follows by the standard semigroup argument.
Let us first solve the equation (\ref{q}). For $n\geq 1$, define recursively that
\begin{align}
q_n(t,x,y):=\int_0^t\!\!\!\int_{\mR^d}q_0(t-s,x,z)q_{n-1}(s,z,y)\dif z\dif s. \label{qn}
\end{align}
We have:

\bl\label{Le3}
There exists a constant $C_d>0$ such that for all $n\in\mN$,
\begin{align}
|q_n(t,x,y)|\leq\frac{(C_d\Gamma(\beta))^{n+1}}{\Gamma((n+1)\beta)}\left(\varrho^0_{(n+1)\beta}(t,x-y)
+\varrho^\beta_{n\beta}(t,x-y)\right),\label{eqn}
\end{align}
where $\Gamma$ is the usual Gamma function.
\el

\begin{proof}
First of all, by our assumptions and (\ref{frp0}), it is easy to see that
$$
|q_0(t,x,y)|\leq C_d\varrho_0^\beta(t,x-y).
$$
Notice that
$$
\cB(\gamma,\beta) \mbox{ is symmetric and non-increasing with respect to each
variable $\gamma$ and $\beta$}.
$$
For $n=1$, by Lemma \ref{3p}, we deduce that
\begin{align*}
|q_1(t,x,y)|&\leq C_d^2\!\!\int_0^t\!\!\!\int_{\mR^d}\varrho_0^\beta(t-s,x-z)\varrho_0^\beta(s,z-y)\dif z\dif s\\
&\leq C_d^2\cB(2\beta,1)\varrho^0_{2\beta}(t,x-y)+C_d^2\cB(\beta,1)\varrho^\beta_\beta(t,x-y)\\
&\leq C_d^2\cB(\beta,\beta)\Big(\varrho^0_{2\beta}+\varrho^\beta_\beta\Big)(t,x-y).
\end{align*}
Suppose now that
$$
|q_n(t,x,y)|\leq\gamma_n\Big(\varrho^0_{(n+1)\beta}+\varrho^\beta_{n\beta}\Big)(t,x-y),
$$
where $\gamma_n>0$ will be determined below.
Using Lemma \ref{3p} again, we have
\begin{align*}
|q_{n+1}(t,x,y)|&\leq C_d\gamma_n\Big(\cB\big(\beta,1+(n+1)\beta\big)+\cB\big((n+2)\beta,1\big)+\cB(2\beta,1+n\beta)\Big)\varrho^0_{(n+2)\beta}(t,x-y)\\
&\quad+C_d\gamma_{n} \Big(\cB\big((n+1)\beta,1\big)+\cB(\beta,1+n\beta)\Big)\varrho^\beta_{(n+1)\beta}(t,x-y)\\
&\leq C_d\gamma_{n}
\cB\big(\beta,(n+1)\beta\big)\Big(\varrho^0_{(n+2)\beta}+\varrho^\beta_{(n+1)\beta}\Big)(t,x-y)\\
&=:\gamma_{n+1}\Big(\varrho^0_{(n+2)\beta}+\varrho^\beta_{(n+1)\beta}\Big)(t,x-y),
\end{align*}
where
$$
\gamma_{n+1}=C_d\gamma_{n} \cB\big(\beta,(n+1)\beta\big).
$$
Hence, by
$\cB(\gamma,\beta)=\frac{\Gamma(\gamma)\Gamma(\beta)}{\Gamma(\gamma+\beta)}$,
we obtain
\begin{align*}
\gamma_{n}= C_d^{n+1}\cB(\beta,\beta)\cB(\beta,2\beta)\cdots\cB(\beta,n\beta)=
\frac{(C_d\Gamma(\beta))^{n+1}}{\Gamma((n+1)\beta)},
\end{align*}
which gives (\ref{eqn}). The proof is complete.
\end{proof}

We also need the H\"older continuity of $q_n$ with respect to $x$.
\bl\label{qnh}
For all $n\geq 0$ and $\gamma\in(0,\beta)$, there exists a constant $C_d>0$ such that
\begin{align*}
&|q_n(t,x,y)-q_n(t,x',y)|\leq \frac{(C_d\Gamma(\beta))^{n+1}}{\Gamma(n\beta+\gamma)}\Big(|x-x'|^{\beta-\gamma}\wedge 1\Big)\\
&\quad\times\Bigg\{\Big(\varrho^0_{\gamma+n\beta}+\varrho^\beta_{\gamma+(n-1)\beta}\Big)(t,x-y)
+\Big(\varrho^0_{\gamma+n\beta}+\varrho^\beta_{\gamma+(n-1)\beta}\Big)(t,x'-y)\Bigg\}.
\end{align*}
\el
\begin{proof}
Let us first prove the following estimate: for every $\gamma\in(0,\beta)$,
\begin{align}
&|q_0(t,x,y)-q_0(t,x',y)|\no\\
&\qquad\preceq\Big(|x-x'|^{\beta-\gamma}\wedge 1\Big)
\Big\{\big(\varrho^0_\gamma+\varrho^\beta_{\gamma-\beta}\big)(t,x-y)
+\big(\varrho^0_\gamma+\varrho^\beta_{\gamma-\beta}\big)(t,x'-y)\Big\}.\label{eq0}
\end{align}
In the case of $|x-x'|>1$, by (\ref{eqn}) we have
$$
|q_0(t,x,y)|\preceq \varrho^\beta_0(t,x-y)
\leq\varrho^\beta_{\gamma-\beta}(t,x-y)
$$
and
$$
|q_0(t,x',y)|\preceq \varrho^\beta_0(t,x'-y)
\leq\varrho^\beta_{\gamma-\beta}(t,x'-y).
$$
If $t<|x-x'|\leq1$, we can also deduce by (\ref{eqn}) that
$$
|q_0(t,x,y)|\preceq \varrho^\beta_0(t,x-y)
\leq|x-x'|^{\beta-\gamma}\varrho^\beta_{\gamma-\beta}(t,x-y),
$$
and
$$
|q_0(t,x',y)|\preceq|x-x'|^{\beta-\gamma}\varrho^\beta_{\gamma-\beta}(t,x'-y).
$$
Suppose now that
\begin{align*}
|x-x'|\leq t.
\end{align*}
Without loss of generality, we may assume that $x$ is nearer to $y$, i.e.,
$$
|x-y|\leq |x'-y|.
$$
By the definition of $q_0$, we can write
\begin{align*}
|q_0(t,x,y)-q_0(t,x',y)|&\leq |\big(\sL^{\kappa(x)}-\sL^{\kappa(x')}\big)p_0(t,x',y)|\\
&+|\big(\sL^{\kappa(x)}-\sL^{\kappa(y)}\big)\big(p_0(t,x,y)-p_0(t,x',y)\big)|\\
&+|b(x)-b(x')|\cdot|\nabla p_0(t,x',y)|\\
&+|b(x)-b(y)|\cdot|\nabla p_0(t,x,y)-\nabla p_0(t,x',y)|\\
&=:\cI_1+\cI_2+\cI_3+\cI_4.
\end{align*}
Using (\ref{kappa2}) and (\ref{frp0}), we have
$$
\cI_1\preceq |x-x'|^{\beta}\varrho^0_0(t,x'-y)\leq |x-x'|^{\beta-\gamma}\varrho^0_\gamma(t,x'-y).
$$
For the second term, taking $\vartheta=\beta-\gamma$ in (\ref{p03}) yields that
\begin{align*}
\cI_2\preceq (|x-y|^{\beta}\wedge1)\cdot|x-x'|^{\beta-\gamma}\varrho^0_{\gamma-\beta}(t,x-y)=|x-x'|^{\beta-\gamma}\varrho^\beta_{\gamma-\beta}(t,x-y).
\end{align*}
For $\cI_3$, it holds by (\ref{frp0}) that
$$
\cI_3\preceq |x-x'|^{\beta}\varrho^0_{0}(t,x'-y)\leq |x-x'|^{\beta-\gamma}\varrho^0_{\gamma}(t,x'-y).
$$
As for the last term, it follows by taking $\vartheta=\beta-\gamma$ in (\ref{p02}) that
$$
\cI_4\preceq (|x-y|^{\beta}\wedge1)\cdot|x-x'|^{\beta-\gamma}\varrho^0_{\gamma-\beta}(t,x-y)=|x-x'|^{\beta-\gamma}\varrho^\beta_{\gamma-\beta}(t,x-y).
$$
Combining the above calculations, we obtain (\ref{eq0}).

Now, by the definition of $q_n$ and Lemma \ref{Le3}, we have for $n\in\mN$,
\begin{align*}
&|q_n(t,x,y)-q_n(t,x',y)|\preceq\int_0^t\!\!\!\int_{\mR^d}|q_0(t-s,x,z)-q_0(t-s,x',z)|\cdot|q_{n-1}(s,z,y)|\dif z\dif s\\
&\quad\preceq\frac{(C_d\Gamma(\beta))^{n}}{\Gamma(n\beta)}\Big(|x-x'|^{\beta-\gamma}\wedge 1\Big)\int_0^t\!\!\!\int_{\mR^d}
\Big\{\big(\varrho_{\gamma}^0+\varrho_{\gamma-\beta}^{\beta}\big)(t-s,x-z)+\big(\varrho_{\gamma}^0+\varrho_{\gamma-\beta}^{\beta}\big)(t-s,x'-z)\Big\}\\
&\qquad\qquad\qquad\qquad\qquad\qquad\qquad\times\Big\{\varrho_{n\beta}^{0}(s,z-y)+\varrho_{(n-1)\beta}^{\beta}(s,z-y)\Big\}\dif z\dif s,
\end{align*}
which yields the desired result by Lemma \ref{3p}.
\end{proof}

Basing on the above two lemmas, we have
\bl\label{T3.4}
The function $q(t,x,y):=\sum_{n=0}^{\infty}q_n(t,x,y)$ solves
the integro-differential equation (\ref{q}). Moreover, $q(t,x,y)$ has the following estimates:
\begin{align}
|q(t,x,y)|\preceq\varrho^\beta_0(t,x-y)+\varrho^0_\beta(t,x-y), \label{eq3}
\end{align}
and any $\gamma\in(0,\beta)$,
\begin{align}
&|q(t,x,y)-q(t,x',y)|\no\\
&\qquad\preceq\Big(|x-x'|^{\beta-\gamma}\wedge 1\Big)\Big\{(\varrho^0_\gamma+\varrho^\beta_{\gamma-\beta})(t,x-y)
+(\varrho^0_\gamma+\varrho^\beta_{\gamma-\beta})(t,x'-y)\Big\}.
\label{eq4}
\end{align}
\el

\begin{proof}
By Lemma \ref{Le3}, one sees that
\begin{align*}
\sum_{n=0}^{\infty}|q_n(t,x,y)|&\leq\sum_{n=0}^\infty
\frac{(C_d\Gamma(\beta))^{n+1}}{\Gamma((n+1)\beta)}\left(\varrho^0_{(n+1)\beta}(t,x,y)
+\varrho^\beta_{n\beta}(t,x,y)\right)\\
&\leq\left\{\sum_{n=0}^\infty\frac{(C_d\Gamma(\beta))^{n+1}}{\Gamma((n+1)\beta)}\right\}
\left(\varrho^0_{\beta}(t,x,y)+\varrho^\beta_0(t,x,y)\right).
\end{align*}
Since the series is convergent, we obtain (\ref{eq3}). Similarly, estimate (\ref{eq4}) follows by Lemma \ref{qnh}.
Moreover, by (\ref{qn}) we have
$$
\sum_{n=0}^{m+1}q_n(t,x,y)=q_0(t,x,y)+\int^t_0\!\!\!\int_{\mR^d}q_0(t-s,x,z)\sum_{n=0}^mq_n(s,z,y)\dif z\dif s,
$$
which yields (\ref{q}) by taking limits $m\to\infty$ for both sides.
\end{proof}

For brevity, set
$$
\sQ(t,x,y):=\int_0^t\!\!\!\int_{\mR^d}p_0(t-s,x,z)q(s,z,y)\dif z\dif s.
$$
With Lemma \ref{T3.4} in hand, we can prove the following results, whose proof is entirely similar to the one of \cite[Lemma 3.4, Lemma 3.5]{Xi-Zh}, we omit the details here. See also \cite[Lemma 3.5, Lemma 3.6]{Ch-Zh}.

\bl\label{omit}
For all $t>0$ and $x\not=y\in\mR^d$, we have
\begin{align*}
\p_t \sQ(t,x,y)=-q(t,x,y)-\int^t_0\!\!\!\int_{\mR^d}\sL^{\kappa(y)}p_0(t-s,\cdot,z)(x)q(s,z,y)\dif z\dif s,
\end{align*}
and
\begin{align*}
\nabla_x \sQ(t,x,y)&=\int_0^t\!\!\!\int_{\mR^d}\nabla_xp_0(t-s,x,z)q(s,z,y)\dif z\dif s,\\
\sL^{\kappa(x)} \sQ(t,x,y)&=\int^t_0\!\!\!\int_{\mR^d}\sL^{\kappa(x)}p_0(t-s,\cdot,z)(x)q(s,z,y)\dif z\dif s,
\end{align*}
where the integrals are understood in the sense of iterated integrals.
\el

Before giving the proof of the main result, we prepare the following non-local maximal principle, see \cite{Xi-Zh} and \cite{Ch-Zh}.

\bt\label{Max}
(Maximal principle)
For given $T>0$, let $u(t,x)\in C_b([0,T]\times\mR^d)$ be such that for almost all $t\in[0,T]$ and all $x\in\mR^d$,
\begin{align}
\p_t u(t,x)+\sL u(t,x)=0.\label{EQQ1}
\end{align}
Assume that
\begin{align}
\lim_{t\uparrow T}\|u(t)-u(T)\|_\infty=0,\ \ \sup_{t\in[0,s]}\|\nabla u(t)\|_\infty<+\infty,\ \ s\in[0,T),\label{EO1}
\end{align}
and
\begin{align}
\mbox{for each $x\in\mR^d$, $t\mapsto \sL u(t,x)$ are continuous on $[0,T)$.}\label{EO11}
\end{align}
Then for each $t\in[0,T)$,
$$
\sup_{x\in\mR^d}u(t,x)\leq\sup_{x\in\mR^d}u(T,x).
$$
In particular, there is a unique solution to equation (\ref{EQQ1}) with the same final value at time $T$ in the class of $u\in C_b([0,T]\times\mR^d)$
satisfying (\ref{EO1}) and (\ref{EO11}).
\et
\begin{proof}
Without loss of generality, we may assume that $u$ is nonnegative. Otherwise, we can subtract  the infimum of $u$ from $u$.
By the assumption, it suffices to prove that for any $t<s<T$,
\begin{align}
\sup_{x\in\mR^d}u(t,x)\leq\sup_{x\in\mR^d}u(s,x).\label{ER911}
\end{align}
Below we fix $s\in(0,T)$.
Let $\chi(x):\mR^d\to[0,1]$ be a smooth function with $\chi(x)=1$ for $|x|\leq 1$ and $\chi(x)=0$ for $|x|>2$. For $R>0$, define the following cutoff function
$$
\chi_R(x):=\chi(x/R).
$$
For $R,\delta>0$, consider
$$
u^\delta_R(t,x):=u(t,x)\chi_R(x)+(t-s)\delta.
$$
Then
\begin{align}
\p_t u^\delta_R(t,x)+\sL u^\delta_R(t,x)=g^\delta_R(t,x)+\delta,\label{EQ66}
\end{align}
where
\begin{align}
g^\delta_R(t,x)&:=\sL^\kappa(u\chi_R)(t,x)-\sL^\kappa u(t,x)\cdot\chi_R(x)+b(x)\cdot \nabla\chi_R(x) u(t,x).\label{UO2}
\end{align}
We proceed to show that for each $\delta>0$, there exists an $R_0\geq 1$ such that for all $t\in[0,s)$ and $R>R_0$,
\begin{align}
\sup_{x\in\mR^d}u^\delta_R(t,x) \leq \sup_{x\in\mR^d}u^\delta_R(s,x). \label{EY1}
\end{align}
If this is proven, then taking $R\to\infty$ and $\delta\to 0$ and noticing that $\sup_{x\in\mR^d}u^\delta_R(s,x)\leq \sup_{x\in\mR^d}u(s,x)$,
we obtain (\ref{ER911}).

We first prove the following claim: For each $s<T$, there exists a constant $C_s>0$ such that
\begin{align}
\sup_{t\in[0,s]}\|g^\delta_R(t)\|_\infty\leq\frac{C_s}{R^{1/2}}.\label{UO1}
\end{align}
{\it Proof of Claim:}
By definition, we have
\begin{align*}
|\sL^\kappa(u\chi_R)(t,x)-\sL^\kappa u(t,x)\cdot\chi_R(x)|&\preceq
\int_{\mR^d}|u(t,x+z)-u(t,x)|~|\chi_R(x+z)-\chi_R(x)|\frac{\dif z}{|z|^{d+1}}\\
&\qquad+|u(t,x)|\int_{\mR^d}|\delta_{\chi_R}(x;z)|\cdot|z|^{-d-1}\dif z\\
&\preceq \|u(t)\|_\infty\|\chi_R\|_\infty^{\frac{1}{2}}\|\nabla\chi_R\|^{\frac{1}{2}}_\infty\int_{|z|>1}\frac{\dif z}{|z|^{d+1/2}}\\
&+\|\nabla u(t)\|_\infty\|\nabla\chi_R\|_\infty\int_{|z|\leq 1}\frac{\dif z}{|z|^{d-1}}+\|u(t)\|_\infty\|\nabla\chi_R\|_\infty\\
&\preceq \|u(t)\|_\infty\|\chi\|_\infty^{\frac{1}{2}}\frac{\|\nabla\chi\|^{\frac{1}{2}}_\infty}{R^{1/2}}
+\|\nabla u(t)\|_\infty\frac{\|\nabla\chi\|_\infty}{R}+\|u(t)\|_\infty\frac{\|\nabla\chi\|_\infty}{R},
\end{align*}
which then gives (\ref{UO1}) by (\ref{UO2}), (\ref{EO1}).

\vspace{2mm}

We now use the contradiction argument to prove (\ref{EY1}). Fix
\begin{align}
R>(C_s/\delta)^2.\label{ER95}
\end{align}
Suppose that (\ref{EY1}) does not hold, since $t\mapsto\sup_{x\in\mR^d} u^\delta_R(t,x)$ is continuous on $[0,s]$,
there must exist a point $t_0\in[0,s)$ such that
$$
\sup_{(t,x)\in[0,s)\times\mR^d}u^\delta_R(t,x)=\sup_{t\in[0,s)}\left(\sup_{x\in\mR^d}u^\delta_R(t,x)\right)=\sup_{x\in\mR^d} u^\delta_R(t_0,x)
$$
and further for some $x_0\in\mR^d$,
$$
\sup_{(t,x)\in[0,s)\times\mR^d}u^\delta_R(t,x)=\sup_{x\in\mR^d} u^\delta_R(t_0,x)=u^\delta_R(t_0,x_0).
$$
In particular,
\begin{align}
\nabla u^\delta_R(t_0,x_0)=0,\label{UO4}
\end{align}
and
\begin{align}
\sL^\kappa u^\delta_R(t_0,x_0)=\lim_{\eps\downarrow 0}\int_{|z|\geq\eps}(u^\delta_R(t_0,x_0+z)-u^\delta_R(t_0,x_0))\kappa(x,z)|z|^{-d-1}\dif z\leq 0.\label{ER97}
\end{align}
Moreover, by (\ref{EQ66}),  for any $h\in(0,s-t_0)$, we have
$$
0\geq \frac{u^\delta_R(t_0+h,x_0)-u^\delta_R(t_0,x_0)}{h}=-\frac{1}{h}\int^{t_0+h}_{t_0}\sL  u^\delta_R(r,x_0)\dif r
+\frac{1}{h}\int^{t_0+h}_{t_0}g^\delta_R(r,x_0)\dif r+\delta.
$$
Since
$$
t\mapsto \sL u^\delta_R(t,x_0)\mbox{ are continuous},
$$
letting $h\to 0$,  by (\ref{UO4}), (\ref{ER97}) and (\ref{UO1}), we obtain
$$
0\geq\sL u^\delta_R(t_0,x_0)-\frac{C_s}{R^{1/2}}+\delta\geq
-\frac{C_s}{R^{1/2}}+\delta,
$$
which produces a contradiction with (\ref{ER95}). The proof is complete.
\end{proof}

Now, we prove the first main result of this paper.

\begin{proof}[Proof of Theorem \ref{main1}]

First of all, by Lemma \ref{omit}, one sees that the computations in (\ref{form}) make sense, and thus (\ref{eq15}) is true. Meanwhile, the uniqueness, non-negative, conservativeness and the semigroup properties can be obtained by Theorem \ref{Max} with the same arguments as in \cite{Xi-Zh}. Thus, $p(t,x,y)$ forms a density function. We only need to prove the corresponding estimates.\\

\noindent(i) Recalling that $t\in(0,1)$, one has by (\ref{eq3}) and (\ref{3p})
\begin{align*}
|\sQ(t,x,y)|
&\preceq \int_0^t\!\!\!\int_{\mR^d}\varrho_1^0(t-s,x,z)\Big(\varrho^\beta_0+\varrho^0_\beta\Big)(s,z-y)\dif z\dif s\\
&\preceq\varrho^0_{1+\beta}(t,x-y)+\varrho^{\beta}_{1}(t,x-y)\leq\varrho^0_1(t,x-y),
\end{align*}
which in turn gives estimate (\ref{eq16}) by equation (\ref{defp}) and (\ref{p00}).\\

\noindent(ii) We write
\begin{align*}
\nabla \sQ(t,x,y)&=\int^t_{\frac{t}{2}}\!\!\!\int_{\mR^d}\nabla p_0(t-s,x,z)\Big(q(s,z,y)-q(s,x,y)\Big)\dif z\dif s\\
&\quad+\int^t_{\frac{t}{2}}\left(\int_{\mR^d}\nabla p_0(t-s,x,z)\dif z\right)q(s,x,y)\dif s\\
&\quad+\int_0^{\frac{t}{2}}\!\int_{\mR^d}\nabla p_0(t-s,x,z)q(s,z,y)\dif z\dif s\\
&=:\sQ_1(t,x,y)+\sQ_2(t,x,y)+\sQ_3(t,x,y).
\end{align*}
For $\sQ_1(t,x,y)$, by (\ref{frp0}), (\ref{eq4}) and Lemma \ref{3p}, we have
\begin{align*}
|\sQ_1(t,x,y)|&\preceq\int^t_{\frac{t}{2}}\!\int_{\mR^d}
\varrho^{\beta-\gamma}_0(t-s,x-z)\Big\{(\varrho^0_\gamma+\varrho^\beta_{\gamma-\beta})(s,x-y)
+(\varrho^0_\gamma+\varrho^\beta_{\gamma-\beta})(s,z-y)\Big\}\dif z\dif s\\
&\leq\int^t_{\frac{t}{2}}\left(\int_{\mR^d}\varrho^{\beta-\gamma}_0(t-s,x-z)\dif z\right)
(\varrho^0_\gamma+\varrho^\beta_{\gamma-\beta})(s,x-y)\dif s\\
&\quad+\int^t_0\!\!\!\int_{\mR^d}\varrho^{\beta-\gamma}_0(t-s,x-z)(\varrho^0_\gamma+\varrho^\beta_{\gamma-\beta})(s,z-y)\dif z\dif s\\
&\preceq\left(\int^t_{\frac{t}{2}}(t-s)^{\beta-\gamma-1}
(1+s^{\gamma-\beta})\varrho^0_0(s,x-y)\dif s\right)\\
&\quad+(\varrho^0_\beta+\varrho^\beta_{0}+\varrho^{\beta-\gamma}_{\gamma})(t,x-y)\preceq\varrho^0_0(t,x-y).
\end{align*}
Thanks to (\ref{00}),  we can deduce for the second term that
\begin{align*}
|\sQ_2(t,x,y)|\stackrel{(\ref{eq3})}{\preceq}\int^t_{\frac{t}{2}}(t-s)^{\beta-1}
\Big\{\varrho^\beta_0(s,x-y)+\varrho^0_\beta(s,x-y)\Big\}\dif s
\preceq\varrho^0_0(t,x-y).
\end{align*}
As for $\sQ_3(t,x,y)$, we have
\begin{align*}
|\sQ_3(t,x,y)|&\preceq\int_0^{\frac{t}{2}}\!\int_{\mR^d}
\varrho^0_0(t-s,x-z)\Big\{\varrho^\beta_0(s,z-y)+\varrho^0_\beta(s,z-y)\Big\}\dif z\dif s\\
&\preceq t^{-1}\!\!\int_0^{t}\!\int_{\mR^d}
\varrho^0_1(t-s,x-z)\Big\{\varrho^\beta_0(s,z-y)+\varrho^0_\beta(s,z-y)\Big\}\dif z\dif s
\preceq\varrho^0_0(t,x-y).
\end{align*}
Combining the above calculations, we obtain
\begin{align*}
|\nabla \sQ(t,x,y)|\preceq\varrho^0_0(t,x-y),
\end{align*}
which in turn gives (\ref{eq17}).

\noindent(iii) Set
$$
\zeta_{\nabla p_0}(t;x,x';y):=\nabla p_0(t,x,y)-\nabla p_0(t,x',y).
$$
Then, estimate (\ref{p02}) yields that for any $\vartheta\in(0,1)$,
$$
|\zeta_{\nabla p_0}(t;x,x';y)|\leq C_\vartheta|x-x'|^\vartheta\varrho_{-\vartheta}^0(t,\tilde x-y),
$$
where $\tilde x$ is the one of the two points $x$ and $x'$ which is nearer to $y$.
We may argue as above to write
\begin{align*}
\big|\nabla\sQ(t,x,y)-\nabla\sQ(t,x',y)\big|&\leq \Bigg|\int_{\frac{t}{2}}^t\!\!\!\int_{\mR^d}\zeta_{\nabla p_0}(t-s;x,x';z)\Big(q(s,z,y)-q(s,\tilde x,y)\Big)\dif z\dif s\Bigg|\\
&\quad+\Bigg|\int_{\frac{t}{2}}^t\!\!\!\int_{\mR^d}\zeta_{\nabla p_0}(t-s;x,x';z)\dif zq(s,\tilde x,y)\dif s\Bigg|\\
&\quad+\Bigg|\!\int^{\frac{t}{2}}_0\!\!\!\int_{\mR^d}\zeta_{\nabla p_0}(t-s;x,x';z)q(s,z,y)\dif z\dif s\Bigg|\\
&=:\sD_1(t,x,x',y)+\sD_2(t,x,x',y)+\sD_3(t,x,x',y).
\end{align*}
For $\vartheta<\beta$, we can choose a $\vartheta'>0$ such that $\vartheta+\vartheta'<\beta$, and by (\ref{eq4}), (\ref{3p}) we have
\begin{align*}
\sD_1(t,x,x',y)&\leq C_\vartheta|x-x'|^\vartheta\!\!\int_{\frac{t}{2}}^t\!\!\!\int_{\mR^d}\varrho^{\beta-\vartheta'}_{-\vartheta}(t-s,\tilde x-z)\Big(\varrho^0_{\vartheta'}+\varrho^{\beta}_{\vartheta'-\beta}\Big)(s,z-y)\dif z\dif s\\
&\quad+C_\vartheta|x-x'|^\vartheta\!\!\int_{\frac{t}{2}}^t\!\!\!\int_{\mR^d}\varrho^{\beta-\vartheta'}_{-\vartheta}(t-s,\tilde x-z)\dif z\Big(\varrho^0_{\vartheta'}+\varrho^{\beta}_{\vartheta'-\beta}\Big)(s,x-y)\dif s\\
&\preceq |x-x'|^\vartheta\Big(\varrho^{0}_{\beta-\vartheta}+\varrho^{\beta-\vartheta'}_{\vartheta'-\vartheta}+\varrho^{\beta}_{-\vartheta}\Big)(t,\tilde x-y)\\
&\quad+|x-x'|^\vartheta\!\!\int_{\frac{t}{2}}^t(t-s)^{\beta-\vartheta-\vartheta'-1}\Big(\varrho^0_{\vartheta'}+\varrho^{\beta}_{\vartheta'-\beta}\Big)(s,x-y)\dif s\\
&\leq C_\vartheta|x-x'|^\vartheta\varrho_{-\vartheta}^0(t,\tilde x-y).
\end{align*}
Thanks to (\ref{000}) and taken into account of (\ref{eq3}), it holds
\begin{align*}
\sD_2(t,x,x',y)\leq C_\vartheta|x-x'|^\vartheta\!\!\int^t_{\frac{t}{2}}(t-s)^{\beta-\vartheta-1}\Big(\varrho_0^{\beta}+\varrho_\beta^0\Big)(s,\tilde x-y)\dif s\leq C_\vartheta|x-x'|^\vartheta\varrho^0_{-\vartheta}(t,\tilde x-y).
\end{align*}
Finally, we have by (\ref{p02}), (\ref{eq3}) and (\ref{3p}) that for any $\vartheta{'\!'}>0$,
\begin{align*}
\sD_3(t,x,x',y)&\leq C_\vartheta|x-x'|^\vartheta\!\!\int^{\frac{t}{2}}_0\!\!\!\int_{\mR^d}\varrho^0_{-\vartheta}(t-s,x-z)\Big(\varrho_0^{\beta}+\varrho_\beta^0\Big)(s,z-y)\dif z\dif s\\
&\leq C_\vartheta|x-x'|^\vartheta t^{-\vartheta-\vartheta{'\!'}}\!\!\int^{t}_0\!\!\!\int_{\mR^d}\varrho^0_{\vartheta{'\!'}}(t-s,x-z)\Big(\varrho_0^{\beta}+\varrho_\beta^0\Big)(s,z-y)\dif z\dif s\\
&\leq C_\vartheta|x-x'|^\vartheta\varrho^0_{-\vartheta}(t,\tilde x-y).
\end{align*}
Based on the above estimates, we thus get (\ref{f2}) by (\ref{frp0}) and (\ref{defp}). The proof is finished.
\end{proof}

\subsection{Regularity of the semigroup}

At the end of this section, let us consider the following elliptic integral-differential equation in $\mR^d$:
\begin{align}
\lambda u(x)-\sL^\kappa u(x)-b(x)\cdot\nabla u(x)=b(x), \label{pide1}
\end{align}
where $\lambda>0$ is a constant. Denote by $\cT_t$ the semigroup corresponding to $\sL$, i.e.,
$$
\cT_tf(x):=\int_{\mR^d}p(t,x,y)f(y)\dif y,\quad\forall f\in\cB_b(\mR^d).
$$
Using the conclusions obtained above, we can prove the following result.

\bt\label{tt}
Suppose that (\ref{kappa})-(\ref{kappa2}) hold and $b\in C^\beta_b(\mR^d)$ for some $\beta>0$. Then, there exists a classical solution $u\in C_b^{1+\vartheta}(\mR^d)$ to (\ref{pide1}) with $0<\vartheta<\beta$, which is given by
\begin{align*}
u(x)=\int_0^\infty\!\e^{-\lambda t}\cT_tb(x)\dif t.
\end{align*}
Moreover, for $\lambda$ big enough, we have
\begin{align}
\|u\|_{\infty}+\|\nabla u\|_{\infty}\leq \frac{1}{2}.    \label{es2}
\end{align}
\et
\begin{proof}
Recall that
$$
\sL=\sL^\kappa+b\cdot\nabla.
$$
By Fubini's theorem and integral by part formula, we have
\begin{align*}
\sL u(x)&=\int_0^\infty\!\e^{-\lambda t}\sL\cT_tb(x)\dif t=\int_0^\infty\!\e^{-\lambda t}\p_t\cT_tb(x)\dif t\\
&=-b(x)+\lambda u(x),
\end{align*}
which gives (\ref{pide1}). We show that $u\in C_b^{1+\vartheta}(\mR^d)$. As a direct result of (\ref{eq16}), we have
\begin{align}
\|u\|_\infty\leq C_1\lambda^{-1}\|b\|_\infty.  \label{22}
\end{align}
Since $p(t,x,y)$ is a density function, we have
\begin{align}
\int_{\mR^d}p(t,x,y)\dif y=1,\quad\forall x\in\mR^d.  \label{ppp}
\end{align}
As a result, we can write
$$
\nabla u(x)=\int_0^\infty\!\!\!\!\int_{\mR^d}\e^{-\lambda t}\nabla p(t,x,y)\Big(b(y)-b(x)\Big)\dif y\dif t.
$$
Thus, we arrive at
\begin{align*}
\|\nabla u\|_\infty\leq \|b\|_{C^\beta_b}\!\int_0^\infty\!\!\!\!\int_{\mR^d}\e^{-\lambda t}\varrho_0^\beta(t,x-y)\dif y\dif t\leq C\lambda^{-\beta}\|b\|_{C^\beta_b},
\end{align*}
which together with (\ref{22}) implies (\ref{es2}) is true. Finally, using (\ref{ppp}) once more we can write
\begin{align*}
\nabla\cT_tb(x)-\nabla\cT_tb(x')&=\int_{\mR^d}\Big(\nabla p(t,x,y)-\nabla p(t,x',y)\Big)\big(b(y)-b(\tilde x)\big)\dif y,
\end{align*}
where $\tilde x$ is the one of the two points $x$ and $x'$ which is nearer to $y$. In view of (\ref{f2}), we deduce that for $0<\vartheta<\beta$,
\begin{align*}
\nabla\cT_tb(x)-\nabla \cT_tb(x')&\leq C_2|x-x'|^{\vartheta}\|b\|_{C^\beta_b}\!\int_{\mR^d}\varrho_{-\vartheta}^\beta(t,\tilde x-y)\dif y\no\\
&\leq C_2|x-x'|^{\vartheta}t^{\beta-\vartheta-1}\|b\|_{C^\beta_b}.
\end{align*}
Consequently, we find that
\begin{align*}
|\nabla u(x)-\nabla u(x')|&\leq \int_0^\infty\!\e^{-\lambda t}\big|\nabla\cT_tb(x)-\nabla\cT_tb(x')\big|\dif t\\
&\leq C|x-x'|^{\vartheta}\|b\|_{C^\beta_b}\!\int_0^\infty\!t^{\beta-\vartheta-1}\e^{-\lambda t}\dif t\leq C_{\lambda,\beta}|x-x'|^{\vartheta}\|b\|_{C^\beta_b},
\end{align*}
which in turn yields the desired result.
\end{proof}

\section{SDEs driven by Markov process}

In this section, we consider SDE (\ref{sde1}), whose generator $\sL$ is given by (\ref{LLLL}). We want to show the existence and uniqueness of the strong solution of SDE (\ref{sde1}) with irregular coefficients by using the fundamental solution method. Below, we always assume that $\sigma$ satisfies (\ref{kappa})-(\ref{kappa2}) holds with $\beta>\tfrac{1}{2}$, and $b\in C_b^\theta(\mR^d)$ with $\theta>\tfrac{1}{2}$.

\subsection{Krylov estimate and Zvonkin's transformation}
Let us first introduce the following class of functions to be used.
\bd\label{Def1}
(Generalized Kato's class) Define
$$
\mK^1_d:=\Big\{f\in L^1_{\mathrm{loc}}(\mR^d): \text{for every}\,\,T>0, K^1(T)<\infty\Big\},
$$
where
$$
K^1(T):=\sup_{x\in\mR^d} \int_{\mR^d}|f(x-y)| \cdot\left(\frac{1}{|y|^{d-1}}\wedge\frac{T^2}{|y|^{d+1}}\right)\dif y.
$$
\ed
For the characterization for $\mK^1_d$, see \cite{Bo-Ja0} and \cite[Proposition 2.3]{Xi-Zh} for more discussions.
By H\"older's inequality, one can easily see that for $p>d$,
\begin{align*}
K^1(T)\leq \|f\|_p\cdot\!\left(\int_{\mR^d}\bigg(\frac{1}{|y|^{d-1}}\wedge\frac{T^2}{|y|^{d+1}}\bigg)^q\dif y\right)^{\frac{1}{q}}\leq C_T\cdot\!\left(\int_{|y|\leq T}\frac{1}{|y|^{q(d-1)}}\dif y\right)^{\frac{1}{q}}<\infty,
\end{align*}
where $q$ is the conjugate index of $p$ and since $p>d$, we have $q(d-1)<d$. Thus, we get
\begin{align}
L^p(\mR^d)\subseteq \mK^1_d,\quad\forall p>d.    \label{incl}
\end{align}

It was shown in \cite[Proposition 3]{M-P} that under our conditions, there exists a unique martingale solution corresponding to the operator $\sL$. Meanwhile, it is known that the martingale solution for $\sL$ is equivalent to the weak solution to SDE (\ref{sde1}), see \cite[Lemma 2.1]{Kurz2}. Thus, the existence and uniqueness of weak solution hold for SDE (\ref{sde1}). As an application of Theorem \ref{main1}, we have the following result.

\bl
The unique weak solution $X$ of SDE (\ref{sde1}) has a jointly continuous density function $p(t,x,y)$ with respect to the Lebesgue measure on $\mR^d$. Moreover, $p(t,x,y)$ enjoys all the properties stated in Theorem \ref{main1} and for every $T>0$ and any nonnegative function $f\in \mK_d^1$,
\begin{align}
\sup_{x\in\mR^d}\mE\left(\int_0^T\!\!f\big(X_s(x)\big)\dif s\right)\leq C_{d,T}\|f\|_p,    \label{kry2}
\end{align}
where $C_{d,T}$ is a positive constant.
\el

\begin{proof}
The first part of the conclusions follows by the same method as in \cite[Corollary 1.3]{Ch-Zh}. We proceed to show the estimate (\ref{kry2}).
By (\ref{eq16}), we have
\begin{align*}
\mE\left(\int_0^t\!\!f\big(X_s(x)\big)\dif s\right)&=\int_0^t\!\!\!\int_{\mR^d}p(s,x,y)f(y)\dif y\dif s\\
&\leq \int_0^t\!\!\!\int_{\mR^d}\varrho_1^0(s,x-y)f(y)\dif y\dif s\leq \cI_1(t)+\cI_2(t),
\end{align*}
where
$$
\cI_1(t):=\int_0^t\!\!\!\int_{|y|\leq t}\varrho_1^0(s,y)f(x-y)\dif y\dif s,
$$
and
$$
\cI_2(t):=\int_0^t\!\!\!\int_{|y|> t}\varrho_1^0(s,y)f(x-y)\dif y\dif s.
$$
Using the definition of $\mK_d^1$, we find that for $\cI_1(t)$,
\begin{align*}
\cI_1(t)&\leq \int_{|y|\leq t}\left(\int_0^{|y|}\frac{s}{|y|^{d+1}}\dif s+\int_{|y|}^ts^{-d}\dif s\right)f(x-y)\dif y\\
&\leq \int_{|y|\leq t}\frac{1}{|y|^{d-1}}f(x-y)\dif y<\infty.
\end{align*}
As for the second term, we can deduce
\begin{align*}
\cI_2(t)&\leq \int_{|y|> t}\frac{t^2}{|y|^{d+1}}f(x-y)\dif y<\infty.
\end{align*}
The proof is finished.
\end{proof}

\br
Estimate (\ref{kry2}) is called the Krylov estimate for the strong solutions, which is very important and usually obtained by suitable analytic regularity method, see \cite{Kr-Ro,Zh3,Zh00}. Here, we obtain this result by simply using the estimate of the fundamental solution.
\er

Usually, the It\^o's formula is performed for functions $f\in C^2_b(\mR^d)$. However, this is too strong for our latter use. Notice that $\sL^{\kappa}f$ is meaningful for any $f\in C^{1+\gamma}_b(\mR^d)$ as long as $\gamma>0$. Indeed, we have by (\ref{kappa}) that
\begin{align*}
\sL^\kappa f(x)&\leq C_{d}\!\int_{|z|\leq 1}\!\!\int_0^1\!\big|\nabla f(x+rz)-\nabla f(x)\big|\dif r\frac{\dif z}{|z|^{d}}+C_{d}\|f\|_{\infty}\\
&\leq C_{d}\!\int_{|z|\leq 1}\frac{\dif z}{|z|^{d-\gamma}}\|f\|_{1+\gamma}+C_{d}\|f\|_{\infty}<\infty.
\end{align*}
We first show that It\^o's formula holds for $f(X_t)$ when $f\in C^{1+\gamma}_b(\mR^d)$ with $\gamma>0$.

\bl\label{ito}
Let $X_t$ satisfies (\ref{sde1}) and $f\in C^{1+\gamma}_b(\mR^d)$ with $\gamma>0$. Then, we have
\begin{align*}
f(X_t)-f(x)-\int_0^t\!\sL f(X_s)\dif s=\int_0^t\!\!\!\int_0^{\infty}\!\!\!\!\int_{\mR^d}\big[f\big(X_{s-}+1_{[0,\sigma(X_{s-},z)]}(r)z\big)-f(X_s)\big]\tilde N(\dif z\times\dif r\times \dif s).
\end{align*}
\el
\begin{proof}
Let $\rho\in C^{\infty}_0(\mR^d)$ such that $\int_{\mR^d}\rho(x)\dif x=1$. Define $\rho_n(x):=n^d\rho(nx)$, and
\begin{align*}
f_n(x):=\int_{\mR^d}f(y)\rho_n(x-y)\dif y.
\end{align*}
Hence, we have $f_n\in C^2_b(\mR^d)$ with $\|f_n\|_{C^\gamma_b}\leq \|f\|_{C^\gamma_b}$, and $\|f_n-f\|_{C^{\gamma'}_b}\rightarrow0$ for every $\gamma'<\gamma$. By using It\^o's formula for $f_n(X_t)$, we get
\begin{align*}
f_n(X_t)-f_n(x)-\int_0^t\!\sL f_n(X_s)\dif s=\int_0^t\!\!\!\int_0^{\infty}\!\!\!\!\int_{\mR^d}\big[f_n\big(X_{s-}+1_{[0,\sigma(X_{s-},z)]}(r)z\big)-f_n(X_{s-})\big]\tilde N(\dif z\times\dif r\times \dif s).
\end{align*}
Now we are going to pass the limits on the both sides of the above equality. It is easy to see that for every $\omega$ and $x\in\mR^d$,
$$
f_n(X_t)-f_n(x)\rightarrow f(X_t)-f(x),\quad \text{as}\,\,n\rightarrow\infty.
$$
Since
\begin{align*}
|f_n(x+z)-f_n(x)-z\cdot\nabla f_n(x)|\leq C|z|^{\gamma}\|f_n\|_{C^\gamma_b}\leq C|z|^{\gamma}\|f\|_{C^\gamma_b},
\end{align*}
we can get by dominated convergence theorem that for every $\omega$,
$$
\int_0^t\!\sL f_n(X_s)\dif s\rightarrow\int_0^t\!\sL f(X_s)\dif s,\quad \text{as}\,\,n\rightarrow\infty.
$$
Finally, by the isometry formula, we have
\begin{align*}
&\mE\bigg|\!\int_0^t\!\!\!\int_0^{\infty}\!\!\!\!\int_{\mR^d}\Big[f_n\big(X_{s-}+1_{[0,\sigma(X_{s-},z)]}(r)z\big)-f_n(X_s)\\ &\qquad-f\big(X_{s-}+1_{[0,\sigma(X_{s-},z)]}(r)z\big)+f(X_{s-})\Big]\tilde N(\dif z\times\dif r\times \dif s)\bigg|^2\\
&=\mE\int_0^t\!\!\!\int_{\mR^d}\!\!\int_0^{\infty}1_{[0,\sigma(X_s,z)]}(r)\big|f_n(X_s+z)-f_n(X_s) -f(X_s+z)+f(X_s)\big|^2\dif r\nu(\dif z)\dif s\\
&\leq C\!\!\int_0^t\!\!\!\int_{\mR^d}\mE\big|f_n(X_s+z)-f_n(X_s) -f(X_s+z)+f(X_s)\big|^2\nu(\dif z)\dif s\rightarrow0,\quad \text{as}\,\,n\rightarrow\infty,
\end{align*}
where in the last step we have used the fact that $\sigma$ is bounded, $\|f_n\|_{C^\gamma_b}\leq \|f\|_{C^\gamma_b}$ and the dominated convergence theorem again. The proof is finished.
\end{proof}

Now, let $u$ be the solution to equation (\ref{pide1}) corresponding to the generator $\sL$ of $X_t$. By Theorem \ref{tt}, we have $u\in C_b^{1+\gamma}(\mR^d)$ with $0<\gamma<\theta\wedge\beta$. Define
$$
\Phi(x):=x+u(x).
$$
In view of (\ref{es2}), we also have
$$
\frac{1}{2}|x-y|\leq\big|\Phi(x)-\Phi(y)\big|\leq \frac{3}{2}|x-y|,
$$
which implies that the map $x\rightarrow\Phi(x)$ forms a $C^1$-diffeomorphism and
\begin{align}
\frac{1}{2}\leq \|\nabla\Phi\|_{\infty},\|\nabla\Phi^{-1}\|_{\infty}\leq 2,   \label{upd}
\end{align}
where $\Phi^{-1}(\cdot)$ is the inverse function of $\Phi(\cdot)$. We prove the following Zvonkin's transformation.

\bl\label{zvon}
Let $\Phi(x)$ be defined as above and $X_t$ solve SDE (\ref{sde1}). Then, $Y_t:=\Phi(X_t)$ satisfies the following SDE:
\begin{align*}
Y_t&=\Phi(x)+\int_0^t\tilde b(Y_s)\dif s+\int_0^t\!\!\!\int_0^{\infty}\!\!\!\!\int_{|z|\leq 1}\tilde g(Y_{s-},z)1_{[0,\tilde\sigma(Y_{s-},z)]}(r)\tilde N(\dif z\times\dif r\times \dif s)\no\\
&\quad+\int_0^t\!\!\!\int_0^{\infty}\!\!\!\!\int_{|z|> 1}\tilde g(Y_{s-},z)1_{[0,\tilde\sigma(Y_{s-},z)]}(r)N(\dif z\times\dif r\times \dif s),
\end{align*}
where
\begin{align*}
\tilde b(x)=\lambda u\big(\Phi^{-1}(x)\big)-\int_{|z|>1}\!\big[u\big(\Phi^{-1}(x)+z\big)-u\big(\Phi^{-1}(x)\big)\big]\sigma\big(\Phi^{-1}(x),z\big)\nu(\dif z)
\end{align*}
and
\begin{align*}
\tilde g(x,z):=\Phi\big(\Phi^{-1}(x)+z\big)-x,\quad \tilde\sigma(x,z):=\sigma\big(\Phi^{-1}(x),z\big).
\end{align*}
\el

\begin{proof}
By Lemma \ref{ito}, we can use the It\^o's formula for function $u$ to get
\begin{align*}
u(X_t)&=u(x)+\int_0^t\!\sL u(X_s)\dif s-\int_0^t\!\!\!\int_0^{\infty}\!\!\!\int_{|z|>1}\big[u\big(X_s+1_{[0,\sigma(X_s,z)]}(v)z\big)-u(X_s)\big]\nu(\dif z)\dif r\dif s\\
&\quad+\int_0^t\!\!\!\int_0^{\infty}\!\!\!\!\int_{|z|> 1}\big[u\big(X_{s-}+1_{[0,\sigma(X_{s-},z)]}(v)z\big)-u(X_{s-})\big] N(\dif z\times\dif r\times \dif s)\\
&\quad+\int_0^t\!\!\!\int_0^{\infty}\!\!\!\!\int_{|z|\leq 1}\big[u\big(X_{s-}+1_{[0,\sigma(X_{s-},z)]}(v)z\big)-u(X_{s-})\big]\tilde N(\dif z\times\dif r\times \dif s).
\end{align*}
Adding this with SDE \eqref{sde1}, taking into account of (\ref{pide1}) and noticing that
$$
\Phi\big(x+y\big)-\Phi(x)=u(x+y)-u(x)+y,
$$
and
\begin{align*}
f\big(x+1_{[0,\sigma(x,z)]}(r)z\big)-f(x)=1_{[0,\sigma(x,z)]}(r)\big[f(x+z)-f(x)\big],
\end{align*}
we obtain the desired result.
\end{proof}

\subsection{Proof of Theorem \ref{main2}}

Before giving the proof of our main results, we prepare some inequalities which will be needed below.

\bl
Let $\tilde b$ and $\tilde g$ be given by Lemma \ref{zvon}. Then, there exist constants $C_1, C_2$ such that for a.e. $x,y\in\mR^d$,
\begin{align}
|\tilde b(x)-\tilde b(y)|\leq C_1|x-y|\cdot\Big(1+h\big(\Phi^{-1}(x)\big)+h\big(\Phi^{-1}(y)\big)\Big)\label{b}
\end{align}
and
\begin{align}
|\tilde g(x,z)-\tilde g(y,z)|\leq C_2|x-y|\cdot|z|^{\gamma},\label{g}
\end{align}
where $0<\gamma<\theta\wedge\beta$.
\el
\begin{proof}
Since $\sigma$ is bounded and thanks to (\ref{es2}), (\ref{upd}), (\ref{a1}), we get
\begin{align*}
|\tilde b(x)-\tilde b(y)|&\leq \lambda\big|u\big(\Phi^{-1}(x)\big)-u\big(\Phi^{-1}(y)\big)\big|+\int_{|z|>1}\!\big|u\big(\Phi^{-1}(x)+z\big)-u\big(\Phi^{-1}(y)+z\big)\big|\nu(\dif z)\\
&\quad+\int_{|z|>1}\!\big|u\big(\Phi^{-1}(x)\big)-u\big(\Phi^{-1}(y)\big)\big|\nu(\dif z)+\int_{|z|>1}\!\big|\sigma\big(\Phi^{-1}(x),z\big)-\sigma\big(\Phi^{-1}(y),z\big)\big|\nu(\dif z)\\
&\leq C_\lambda|x-y|+C_0|x-y|\Big(h\big(\Phi^{-1}(x)\big)+h\big(\Phi^{-1}(y)\big)\Big),
\end{align*}
which gives (\ref{b}). To prove (\ref{g}), we denote by
$$
\cJ_z(x):=u\big(\Phi^{-1}(x)+z\big)-u\big(\Phi^{-1}(x)\big).
$$
Then, one can easily check that
$$
\|\nabla\cJ_z\|_{\infty}\leq C|z|^{\gamma}\|u\|_{C_b^{1+\gamma}}.
$$
Thus, by the definition of $\Phi$ we can deduce
\begin{align*}
|\tilde g(x,z)-\tilde g(y,z)|&=\big|u\big(\Phi^{-1}(x)+z\big)-u\big(\Phi^{-1}(x)\big)-u\big(\Phi^{-1}(y)+z\big)+u\big(\Phi^{-1}(y)\big)\big|\\
&\leq |x-y|\cdot\|\nabla\cJ_z\|_{\infty}\leq C|x-y|\cdot|z|^{\gamma},
\end{align*}
the proof is finished.
\end{proof}

We are now in the position to give:

\begin{proof}[Proof of Theorem \ref{main2}]
Let $X_t$ and $\hat X_t$ be two strong solutions for SDE (\ref{sde1}) both starting from $x\in\mR^d$, and set
$$
Y_t:=\Phi(X_t),\quad \hat Y_t:=\Phi(\hat X_t).
$$
By Lemma \ref{zvon}, we have for all $t\geq 0$,
\begin{align*}
Y_t-\hat Y_t&=\int_0^{t}\!\big[\tilde b(Y_s)-\tilde b(\hat Y_s)\big]\dif s +\!\int_0^{t}\!\!\!\int_0^{\infty}\!\!\!\!\int_{|z|\leq 1}\Big[\tilde g(Y_{s-},z)1_{[0,\tilde\sigma(Y_{s-},z)]}(r)\\
&\quad\qquad\quad\qquad\quad\qquad\quad\quad\qquad-\tilde g(\hat Y_{s-}, z)1_{[0,\tilde\sigma(\hat Y_{s-},z)]}(r)\Big]\tilde N(\dif z\times\dif r\times \dif s)\\
&\quad+\!\int_0^{t}\!\!\!\int_0^{\infty}\!\!\!\!\int_{|z|\leq 1}\Big[\tilde g(Y_{s-},z)1_{[0,\tilde\sigma(Y_{s-},z)]}(r)-\tilde g(\hat Y_{s-}, z)1_{[0,\tilde\sigma(\hat Y_{s-},z)]}(r)\Big] N(\dif z\times\dif r\times \dif s).
\end{align*}
As the argument in \cite[Theorem IV. 9.1]{Ik-Wa} and \cite{Zh00}, we only need to show that
\begin{align}
Z_t\equiv0,   \quad\forall t\geq 0, \label{66}
\end{align}
where $Z_t$ is given by
\begin{align*}
Z_{t}&=\int_0^{t}\!\big[\tilde b(Y_s)-\tilde b(\hat Y_s)\big]\dif s +\!\int_0^{t}\!\!\!\int_0^{\infty}\!\!\!\!\int_{|z|\leq 1}\Big[\tilde g(Y_{s-},z)1_{[0,\tilde\sigma(Y_{s-},z)]}(r)\\
&\quad\qquad\quad\qquad\quad\qquad\quad-\tilde g(\hat Y_{s-}, z)1_{[0,\tilde\sigma(\hat Y_{s-},z)]}(r)\Big]\tilde N(\dif z\times\dif r\times \dif s)=:I^{t}_1+I^{t}_2.
\end{align*}
Set
$$
A(t):=\int_0^t\Big(1+h(X_s)+h(\hat X_s)\Big)\dif s,
$$
then it is easy to see by (\ref{b}) that for and stopping time $\tau$ and almost all $\omega$,
\begin{align*}
\sup_{s\in[0,t]}\left|I^{s\wedge\tau}_1\right|\leq C_1\!\int_0^{t\wedge\tau}|Z_r|\cdot\Big(1+h(X_r)+h(\hat X_r)\Big)\dif r=C_1\!\int_0^{t\wedge\tau}|Z_r| \dif A(r).
\end{align*}
As for the second term, write
\begin{align*}
I_2^{t\wedge\tau}&=\int_0^{t\wedge\tau}\!\!\!\!\int_0^{\infty}\!\!\!\!\int_{|z|\leq 1}1_{[0,\tilde\sigma(Y_{s-},z)\wedge\tilde\sigma(\hat Y_{s-},z)]}(r)\Big[\tilde g(Y_{s-},z)1_{[0,\tilde\sigma(Y_{s-},z)]}(r)\\
&\qquad\qquad\qquad\qquad\qquad\qquad\qquad-\tilde g(\hat Y_{s-}, z)1_{[0,\tilde\sigma(\hat Y_{s-},z)]}(r)\Big]\tilde N(\dif z\times\dif r\times \dif s)\\
&\quad+\int_0^{t\wedge\tau}\!\!\!\!\int_0^{\infty}\!\!\!\!\int_{|z|\leq 1}1_{[\tilde\sigma(Y_{s-},z)\vee\tilde\sigma(\hat Y_{s-},z),\infty]}(r)\Big[\tilde g(Y_{s-},z)1_{[0,\tilde\sigma(Y_{s-},z)]}(r)\\
&\qquad\qquad\qquad\qquad\qquad\qquad\qquad-\tilde g(\hat Y_{s-}, z)1_{[0,\tilde\sigma(\hat Y_{s-},z)]}(r)\Big]\tilde N(\dif z\times\dif r\times \dif s)\\
&\quad+\int_0^{t\wedge\tau}\!\!\!\!\int_0^{\infty}\!\!\!\!\int_{|z|\leq 1}1_{[\tilde\sigma(Y_{s-},z)\wedge\tilde\sigma(\hat Y_{s-},z),\tilde\sigma(Y_{s-},z)\vee\tilde\sigma(\hat Y_{s-},z)]}(r)\Big[\tilde g(Y_{s-},z)1_{[0,\tilde\sigma(Y_{s-},z)]}(r)\\
&\qquad\qquad\qquad\qquad\qquad\qquad\qquad-\tilde g(\hat Y_{s-}, z)1_{[0,\tilde\sigma(\hat Y_{s-},z)]}(r)\Big]\tilde N(\dif z\times\dif r\times \dif s)\\
&=:I_{21}^{t\wedge\tau}+I_{22}^{t\wedge\tau}+I_{23}^{t\wedge\tau}.
\end{align*}
We proceed to estimate each component. First, for $I_{21}^{t\wedge\tau}$, we use the Doob's $L^2$-maximal inequality to deduce that
\begin{align*}
\mE\left[\sup_{s\in[0,t]}|I_{21}^{s\wedge\tau}|\right]&\leq \mE\Bigg(\int_0^{t\wedge\tau}\!\!\!\!\int_0^{\infty}\!\!\!\!\int_{|z|\leq 1}1_{[0,\tilde\sigma(Y_s,z)\wedge\tilde\sigma(\hat Y_s,z)]}(r)\big|\tilde g(Y_s,z)-\tilde g(\hat Y_s, z)\big|^2\dif r\nu(\dif z) \dif s\Bigg)^{\frac{1}{2}}\\
&=\mE\Bigg(\int_0^{t\wedge\tau}\!\!\!\!\int_{|z|\leq 1}\big[\tilde\sigma(Y_s,z)\wedge\tilde\sigma(\hat Y_s,z)\big]\cdot\big|\tilde g(Y_s,z)-\tilde g(\hat Y_s, z)\big|^2\nu(\dif z) \dif s\Bigg)^{\frac{1}{2}}.
\end{align*}
Since $\theta\wedge\beta>1/2$, we can choose $\gamma>0$ such that
$$
1/2<\gamma<\theta\wedge\beta.
$$
We then have by the fact that $\tilde\sigma$ is bounded and (\ref{g}) that
\begin{align*}
\mE\left[\sup_{s\in[0,t]}|I_{21}^{s\wedge\tau}|\right]&\leq C_2\mE\Bigg(\int_0^{t\wedge\tau}|Z_s|^2\!\int_{|z|\leq 1}\!|z|^{2\gamma}\nu(\dif z) \dif s\Bigg)^{\frac{1}{2}}\\
&\leq C_2\mE\Bigg(\int_0^{t\wedge\tau}|Z_s|^2\dif s\Bigg)^{\frac{1}{2}}.
\end{align*}
Next, it is easy to see that for any $t\geq 0$,
$$
I_{22}^{t\wedge\tau}\equiv0.
$$
Finally, we use the $L^1$-estimate (see \cite[P$_{174}$]{Kurz3} or \cite[P$_{157}$]{Kurz2}) to control the third term by
\begin{align*}
\mE\left[\sup_{s\in[0,t]}|I_{23}^{s\wedge\tau}|\right]&\leq 2\mE\!\int_0^{t\wedge\tau}\!\!\!\!\int_0^{\infty}\!\!\!\!\int_{|z|\leq 1}1_{[\tilde\sigma(Y_s,z)\wedge\tilde\sigma(\hat Y_s,z),\tilde\sigma(Y_s,z)\vee\tilde\sigma(\hat Y_s,z)]}(r)\\
&\qquad\qquad\qquad\times\big|\tilde g(Y_s,z)1_{[0,\tilde\sigma(Y_s,z)]}(r)-\tilde g(\hat Y_s, z)1_{[0,\tilde\sigma(\hat Y_s,z)]}(r)\big|\nu(\dif z)\dif r \dif s\\
&\leq2\mE\!\int_0^{t\wedge\tau}\!\!\!\int_{|z|\leq 1}|\tilde\sigma(Y_s,z)-\tilde\sigma(\hat Y_s,z)|\cdot\Big(|\tilde g(Y_s,z)|+|\tilde g(\hat Y_s,z)|\Big)\nu(\dif z)\dif s.
\end{align*}
Since
$$
|\tilde g(x,z)|=\big|\Phi\big(\Phi^{-1}(x)+z\big)-\Phi\big(\Phi^{-1}(x)\big)\big|\leq \frac{3}{2}|z|,
$$
and taking into account of (\ref{a1}), we get
\begin{align*}
\mE\left[\sup_{s\in[0,t]}|I_{23}^{s\wedge\tau}|\right]&\leq C_3\mE\!\int_0^{t\wedge\tau}\!\!\!\int_{|z|\leq 1}|\tilde\sigma(Y_s,z)-\tilde\sigma(\hat Y_s,z)|\cdot|z|\nu(\dif z)\dif s\\
&\leq C_3\mE\!\int_0^{t\wedge\tau}|Z_s|\Big(h(X_s)+h(\hat X_s)\Big)\dif s\leq C_3\mE\!\int_0^{t\wedge\tau_1}|Z_s|\dif A(s).
\end{align*}
Combing the above computations, we arrive at that there exists a constant $C_0$ such that
\begin{align}
\mE\left[\sup_{s\in[0,t]}|Z_{s\wedge\tau}|\right]&\leq C_0\mE\!\int_0^{t\wedge\tau}|Z_s|\dif A(s)+C_0\mE\Bigg(\int_0^{t\wedge\tau}|Z_s|^2\dif s\Bigg)^{\frac{1}{2}}\no\\
&\leq  C_0\mE\!\int_0^{t\wedge\tau}|Z_s|\dif A(s)+C_0\sqrt{t}\cdot\mE\left[\sup_{s\in[0,t]}|Z_{s\wedge\tau}|\right] .\label{ee}
\end{align}
Now, take $t_0$ small enough such that
$$
C_0\sqrt{t_0}<1,
$$
we obtain by (\ref{ee}) that for any stopping time $\tau$,
\begin{align*}
\mE\left[\sup_{s\in[0,t]}|Z_{s\wedge\tau}|\right]&\leq C_1\mE\!\int_0^{t\wedge\tau}\Big[\sup_{r\in[0,s]}|Z_{r}|\Big]\dif A(s).
\end{align*}
By our assumption that $h\in \mK_d^1$ and the Krylov estimate (\ref{kry2}), we find that
$$
\mE A(t)\leq t+C<\infty.
$$
Therefore, $t\mapsto A(t)$ is a continuous strictly increasing process.
As a direct consequence of \cite[Lemma 2.6]{Zh00}, it holds that for almost all $\omega$,
$$
\sup_{s\in[0,t_0]}|Z_{s}|=0.
$$
Since the uniqueness is a local property, we can  get (\ref{66}) by the iteration method. The whole proof is finished.
\end{proof}

\end{document}